\documentclass[11pt]{article}

\usepackage[colorlinks=true]{hyperref}
\usepackage{comment}
\usepackage{amsmath,amssymb}
\usepackage{graphicx,graphics,float}
\usepackage[permil]{overpic}
\usepackage{tikz}
\usepackage{tkz-tab} 
\allowdisplaybreaks
\numberwithin{equation}{section}
\usepackage{amsthm}
\usepackage{enumitem}
\usepackage{esint}

\newtheorem{thm}{Theorem}[section]
\newtheorem{definition}[thm]{Definition}

\newtheorem{remark}[thm]{Remark}
\newtheorem{proposition}[thm]{Proposition}
\newtheorem{corollary}[thm]{Corollary}
\newtheorem{lem}[thm]{Lemma}
\newtheorem{example}[thm]{Example}
\newenvironment{proof*}[1]{\emph{#1}}{\hfill$\blacksquare$}

\usepackage[font=footnotesize, labelfont=normalfont, margin=1in]{caption}
\usepackage[format=plain, margin=1in]{subcaption}
\usepackage{array}
\usepackage{todonotes}

\newcommand{\vt}{\vartheta}

\newcommand{\R}{\mathbb R}
\newcommand{\N}{\mathbb N}
\newcommand{\C}{\mathbb C}

\newcommand{\M}{\mathfrak M}

\renewcommand{\Re}{\mathrm{Re}}

\DeclareMathOperator{\diag}{diag}
\DeclareMathAlphabet{\mathmybb}{U}{bbold}{m}{n}
\newcommand{\1}{\mathmybb{1}}

\usepackage[margin=1in]{geometry}

\definecolor{darkgreen}{HTML}{77AC30}
\definecolor{lightgreen}{RGB}{0 255 0}
\definecolor{myviolet}{HTML}{b240ff}
\definecolor{myblue}{HTML}{0072BD}

\title{Turing instability for nonlocal heterogeneous reaction-diffusion systems:
A computer-assisted proof approach}
\author{Maxime Breden\thanks{CMAP, CNRS, \'Ecole polytechnique, Institut Polytechnique de
Paris, 91120 Palaiseau, France.} $\quad$ Maxime Payan\footnotemark[1] $\quad$ Cordula Reisch\thanks{Institute of Analysis and Scientific Computing, TU Wien, 1040 Vienna, Austria} $\quad$ Bao Quoc Tang\thanks{Department of Mathematics and Scientific Computing, University of Graz, Heinrichstrasse 36, 8010 Graz, Austria}}
\date{\today}
\begin{document}

\maketitle

\abstract{This paper provides a computer-assisted proof for the Turing instability induced by heterogeneous nonlocality in reaction-diffusion systems. Due to the heterogeneity and nonlocality, the linear Fourier analysis gives rise to \textit{strongly coupled} infinite differential systems. 
By introducing suitable changes of basis as well as the Gershgorin disks theorem for infinite matrices, we first show that all $N$-th Gershgorin disks lie completely on the left half-plane for sufficiently large $N$. 
For the remaining finitely many disks, a computer-assisted proof shows that if the intensity $\delta$ of the nonlocal term is large enough, there is precisely one eigenvalue with positive real part, which proves the Turing instability. 
Moreover, by detailed study of this eigenvalue as a function of $\delta$, we obtain a sharp threshold $\delta^*$ which is the bifurcation point for Turing instability. 
}


\tableofcontents

\section{Introduction}

Since Turing's famous paper \cite{turing1952chemical}, the study of 
Turing instability and pattern formation has extensively grown and has become a classical topic in reaction-diffusion systems. Roughly speaking, Turing instability is the phenomenon where spatial diffusion in the short range activator - long range inhibitor regime destabilizes stable homogeneous steady states, leading to the formation of certain patterns, see e.g. \cite{murray2007mathematical} for more details. 
To see this, we usually linearize the system around the homogeneous steady, then by using a Fourier analysis we transform this linearized problem into an infinite \textit{decoupled} differential systems, and finally by studying the spectrum of these systems we can obtain eigenvalues with positive real parts under certain regime of the diffusion, concluding the Turing instability. 
When the system under consideration is heterogeneous or contain non-local terms, this useful procedure might break down, making the study of Turing instability for such systems challenging. 
This also triggers extensive study of Turing instability for heterogeneous systems, see e.g. \cite{banerjee_nonlocal_2022,kozak_pattern_2019,page2003pattern,van_gorder_pattern_2021,chacon2020turing,d2020spatial,banerjee2017spatio}, or non-local systems, see e.g. \cite{colet_formation_2014,taylor2020non,andreguetto_maciel_enhanced_2021,ShiShiSon22,banerjee_nonlocal_2022,pal_nonlocal_2025}. 
To deal with infinite coupling, a heuristic truncation approach is usually employed. 
More precisely, by using the behavior of eigenvalues, it is first shown that the instability of the infinite system can be drawn from the finitely truncated system, provided the truncation is sufficiently large; then numerical simulations are performed on truncated systems to demonstrate and (numerically) confirm the Turing instability. 
As the size of this truncation is usually very large or even not explicitly computable, a theoretical proof for Turing instability for the aforementioned is still out of reach. 
In this paper, we propose the computer-assisted proof to address this issue, and prove the Turing instability for a nonlocal reaction-diffusion system arising from modeling liver inflammation. 
We believe that the proposed approach can be applied to many other situations concerning Turing instability for heterogeneous systems.

\subsection{Problem formulation}
\label{sec:intro_problem}

Let $\Omega$ be a smooth bounded domain of $\mathbb{R}^n$, and $\Omega_1, \Omega_2 \subset \Omega$ be two subdomains.
We consider the following nonlocal reaction-diffusion system
\begin{equation}
\label{eq:system}
\left\{
\begin{aligned}
    &\partial_t u = \vartheta \Delta u + au +bv -u^3, &&x\in\Omega, \\
    &\partial_t v  = \Delta v + cu + dv + \delta \1_{\Omega_1} (x) \frac{1}{|\Omega_1|}\int_{\Omega_2} u \, \mathrm{d}x, &&x\in\Omega,\\
    &\nabla u\cdot \nu = \nabla v \cdot \nu = 0, &&x\in\partial\Omega,
\end{aligned}
\right.
\end{equation}
with parameters $a, b, c, d, \vartheta, \delta \in \mathbb R$, and $\1_{\Omega_1}(\cdot)$ is the characteristic function of $\Omega_1$. 
This model is motivated from modeling liver infection, see e.g. \cite{reisch_chemotactic_2019}, where $u$ and $v$ denote densities of the virus and $\mathsf{T}$~cells, respectively. 
The nonlocal term $\displaystyle \delta \1_{\Omega_1} (x) \frac{1}{|\Omega_1|}\int_{\Omega_2} u$ represents the immune response depending on the virus population in the subdomain $\Omega_2$ through the portal field $\Omega_1$, where $\delta>0$ is the intensity of the response. 
The global existence of solutions to (a variant model of) system \eqref{eq:system} was shown in \cite{reisch_longterm_2022}. 

\begin{figure}
    \centering
    \begin{subfigure}{0.48\textwidth}
            \includegraphics[width=0.98\textwidth]{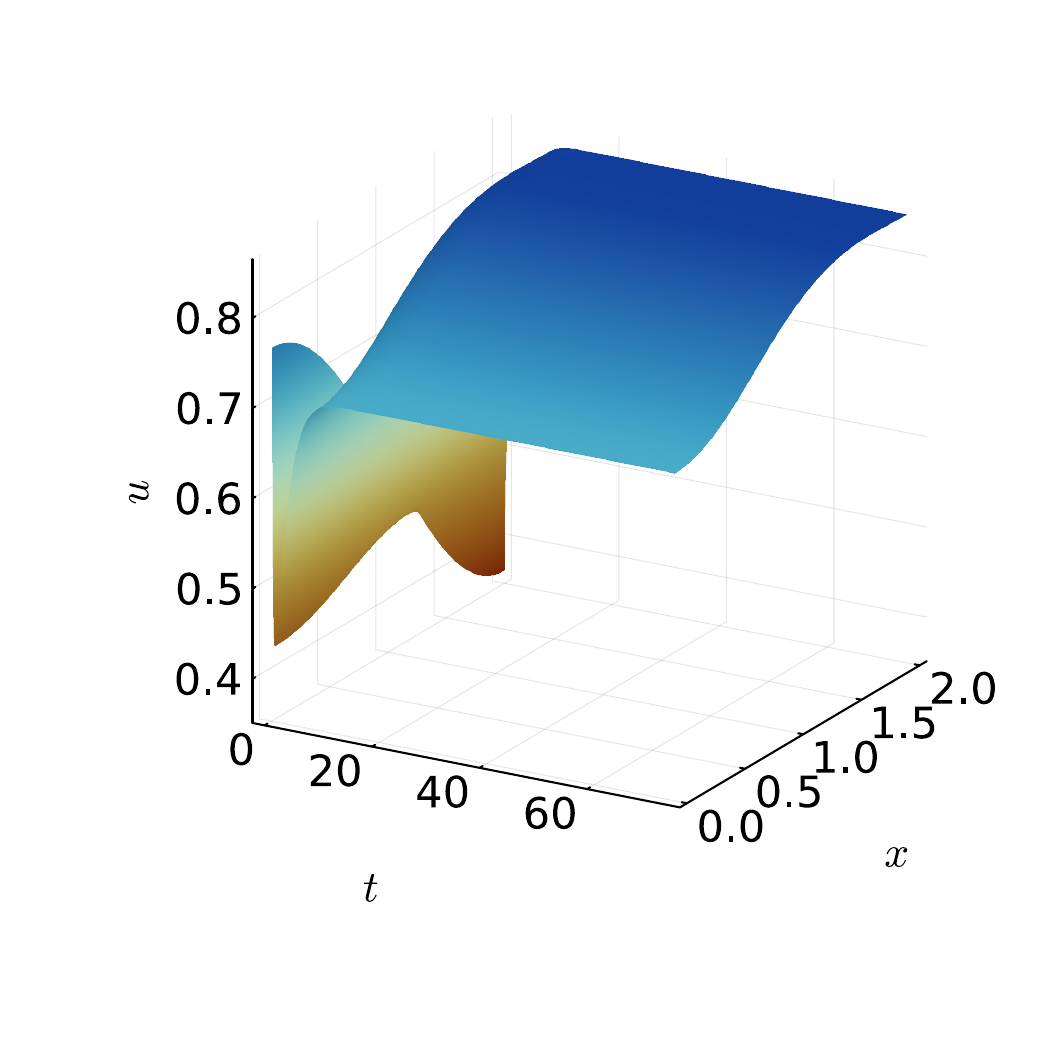}
                \caption{$\delta=4$}
    \label{fig:sol_u_unstable}
    \end{subfigure}
        \begin{subfigure}{0.48\textwidth}
            \includegraphics[width=0.98\textwidth]{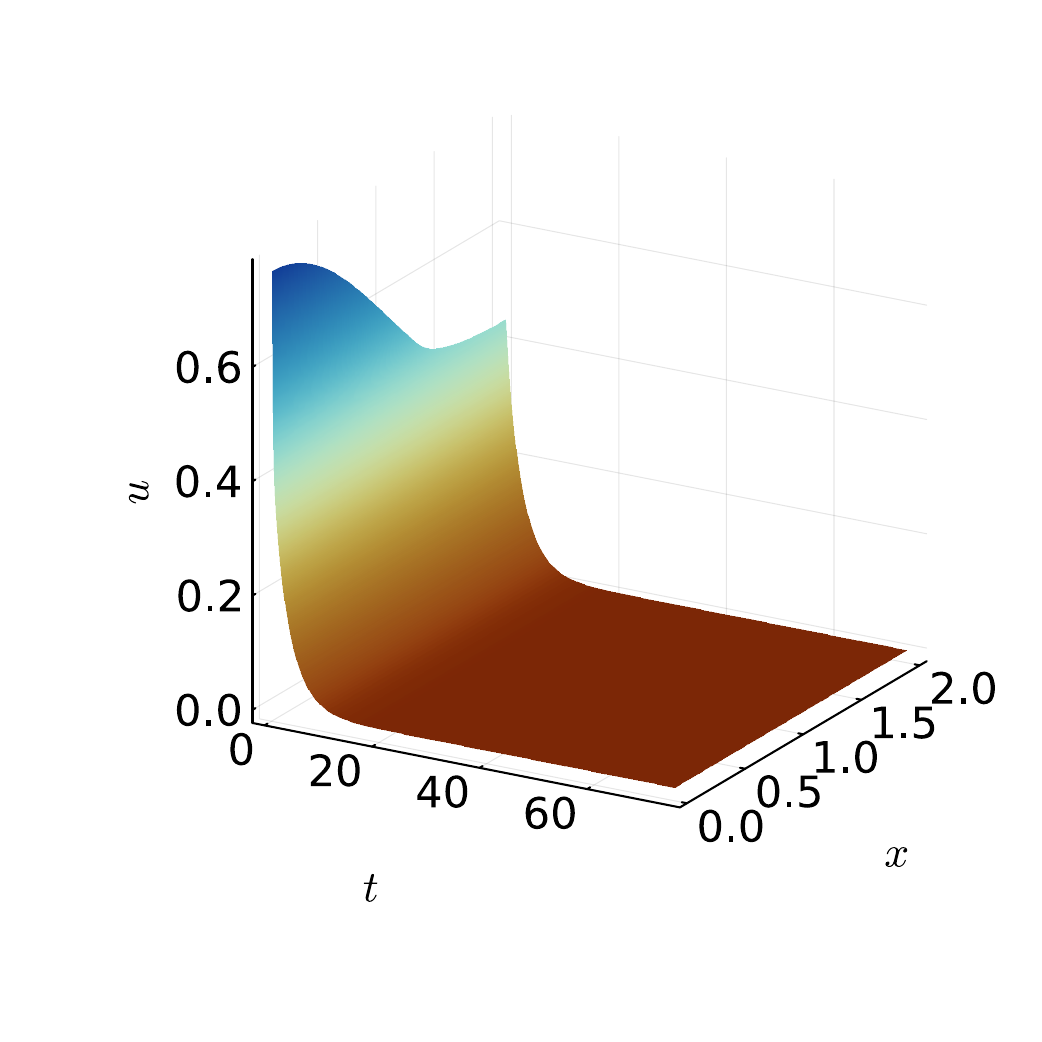}
             \caption{$\delta=1$}
             \label{fig:sol_u_stable}
    \end{subfigure}
    \caption{Numerical simulations of system~\eqref{eq:system}, for the domain and parameter values used in Theorem~\ref{thm:main}. 
    We observe convergence to the trivial equilibrium for $\delta=1$, but the apparition of nontrivial patterns for $\delta=4$, in accordance with Theorem~\ref{thm:main}.}
    \label{fig:solution_u}
\end{figure}

\medskip
System \eqref{eq:system} possesses Turing instability induced by nonlocality. More precisely, let $\delta = 0$ for a moment and consider the parameters $a,b,c,d$ such that the the trivial solution to the corresponding ODE system
\begin{equation*}
    \begin{cases}
        u' = au + bv\\
        v' = cu + dv
    \end{cases}
\end{equation*}
is stable, which requires $a + d < 0$ and $ad - bc > 0$. Moreover, since $u$ and $v$ are densities of virus and $\mathsf{T}$~cells, the solution is expected to be non-negative, and thus it is necessary that $b \geq 0$ and $c\geq 0$. 
Under these conditions, it can be shown, see e.g. \cite{reisch2025Turing}, that the trivial solution to the local PDE system, i.e. \eqref{eq:system} with $\delta = 0$, is also globally stable regardless of the diffusion coefficient $\vartheta>0$. 
This means that the classical Turing instability does not occur. 
Interestingly, it was observed numerically in \cite{reisch2025Turing} that if the intensity $\delta$ is large enough, the trivial solution to \eqref{eq:system} becomes unstable, and the system evolves towards a spatially inhomogeneous state, or in other words, pattern formation appears. 
Fig.~\ref{fig:solution_u} shows pattern formation of \eqref{eq:system} for $\delta = 4$, and in comparison the stable behavior for $\delta=1$. 

To go into more details, we consider the linearized system around the trivial state
\begin{equation}
\label{eq:system_linear1}
\left\{
\begin{aligned}
    &\partial_t u = \vt\Delta u + au +bv, &&x\in\Omega,\\
    &\partial_t v  = \Delta v + cu + dv + \delta \1_{\Omega_1} (x)\frac{1}{|\Omega_1|} \int_{\Omega_2} u \, \mathrm{d}x, &&x\in\Omega,\\
    &\nabla u \cdot \nu = \nabla v \cdot \nu = 0, &&x\in \partial\Omega.
\end{aligned}
\right.
\end{equation}
Let $(\lambda_j,\varphi_j)_{j\geq 0}$ be the eigenvalues-eigenfunctions of the negative Laplacian with homogeneous Neumann boundary conditions, where $0 = \lambda_0 < \lambda_1 \leq \lambda_2 \leq \ldots \to +\infty$ and $\{\varphi_j\}_{j\ge0}$ is an orthonormal basis of $L^2(\Omega)$. 
By writing
\begin{equation*}
    u(x,t) = \sum_{j\geq 0}u_j(t)\varphi_j(x), \quad v(x,t) = \sum_{j\geq 0}v_j(t)\varphi_j(x),
\end{equation*}
it follows from \eqref{eq:system_linear1} that for each $j\geq 0$
\begin{equation}\label{eq:system_infinite}
\left\{
\begin{aligned}
    u_j'
    &= (a- \vartheta\lambda_j) u_j    + b  v_j,\\
    v_j' &= c u_j + (d-\lambda_j) v_j +\delta  \sum_{k\ge0} \left(\frac{1}{|\Omega_1|}\int_{\Omega_1} \varphi_j(x)\mathrm{d}x\int_{\Omega_2} \varphi_k (x)\mathrm{d}x\right)u_k(t).
\end{aligned}
\right.
\end{equation}
It is clear that in general  \eqref{eq:system_infinite} is completely coupled, and therefore obtaining explicitly the solutions seems impossible. 
We denote the infinite vector $X = (X_{j})_{j\geq 0}$ where
\begin{equation*}
    X_{j} = \begin{pmatrix}
        u_{j}\\
        v_{j}
    \end{pmatrix},
\end{equation*}
the infinite matrices
\begin{equation} \label{eq:tildeA}
     \widetilde{A} =     \begin{pmatrix}
        a-\vartheta\lambda_0 & b &  0 & 0 & 0 &  \cdots \\
        c & d- \lambda_0 & 0 & 0 &0 &  \cdots \\
        0 & 0 & a-\vartheta\lambda_1 & b &0 & \cdots \\
          0 & 0 & c & d-\lambda_1 &0 & \cdots\\
        \vdots& \vdots   &\vdots &\vdots & \ddots  &  \vdots\\
    \end{pmatrix}, 
\end{equation}
$B = (B_{i,j})_{i,j\geq 0}$ with
\begin{equation}\label{eq:B}
    B_{i,j} = \frac{1}{|\Omega_1|}\int_{\Omega_1} \varphi_i(x)\mathrm{d}x\int_{\Omega_2} \varphi_j (x)\mathrm{d}x,
\end{equation}
and $\widetilde{B} = (\widetilde{B}_{i,j})_{i,j\geq 0}$ with 
\begin{equation} \label{eq:B1}
    \widetilde{B}_{i,j}= \left\{
    \begin{aligned}
        &B_{\frac{i-1}{2},\frac{j}{2}}, &&\text{if $i$ is odd and $j$ is even}, \\
        &0, &&\text{otherwise.}
    \end{aligned}
    \right.
\end{equation}
Then \eqref{eq:system_infinite} can be written as
\begin{equation}\label{eq:amplitudes_matrix_short}
X' = \widetilde{A}X +\delta \widetilde{B}X =: M X, 
\end{equation}
where $\widetilde{A}$ gives the dynamic of the system without the nonlocal term and $\widetilde{B}$ gives the interaction caused by the nonlocality. 
Investigating the spectrum of the infinite matrix $M$ is a highly non-trivial task. 
Heuristically, since $\lambda_j \underset{j\to+\infty}{\longrightarrow} +\infty$, one can expect that if we truncate the matrix $M$ to $M_N \in \mathbb R^{2N\times 2N}$ for sufficiently large $N$, then when $M_N$ possesses eigenvalues with positive real part would imply the same $M$, leading to the desired instability. 
This was argued in e.g. \cite{kozak_pattern_2019} where the eigenvalues of the truncated matrix were checked numerically. 
In this paper, we propose a rigorous approach to address this issue by Gershgorin disks theorems and computer-assisted proofs.

Before moving on to the heart of the paper, let us mention that there are many other works using computer-assisted proofs to study PDEs and their dynamics. In the specific context of computer-assisted proofs, the Gershgorin theorem has already been used, see, e.g.,~\cite{Rum20,GomOrr21,LesPug24,BreChu25}, but only for finite dimensional matrices, whereas we deal here with the infinite matrix $M$, meaning that truncations errors have to be controlled rigorously.  For a broader overview about computer-assisted proofs in PDEs and other dynamical systems, and many other applications of the techniques used in Section~\ref{sec:NK}, we refer to the expository article~\cite{BerLes15}, the survey~\cite{Gom19}, the book~\cite{NakPluWat19} and the references therein.  


While we focus on system~\eqref{eq:system_linear1} as a specific example (see Theorem~\ref{thm:main}), the general strategy presented in the paper, as well as many intermediate results, remain valid for a wider class of nonlocal terms than $\1_{\Omega_1} (x)\frac{1}{|\Omega_1|} \int_{\Omega_2} u \, \mathrm{d}x$. In particular, the general estimates derived in Section~\ref{sec:Gersh} to control the spectrum of $M$ apply to any nonlocal term of the form $\int_\Omega g(x,y) u(y) \mathrm{d}y$, provided the corresponding operator $B$, now defined by
\begin{equation}
    B_{i,j} = \int_{\Omega}\int_{\Omega} g(x,y)\varphi_i(x)\varphi_j(y)\mathrm{d}x\mathrm{d}y,
\end{equation}
satisfies the following assumption:
\begin{align} \label{hyp:B}
    \exists q_1>\frac{1}{2}, \exists q_2 > -\frac{3}{2} , \exists C\geq 0, \ \forall i,j \geq 0, \ |B_{i,j}| \leq \dfrac{C}{\max(1,i^{q_1})\max(1,j^{q_2})}. \tag{H:\textit{B}}
\end{align}
For examples of nonlocal terms appearing in various models, see~\cite{pal_nonlocal_2025} and the references therein. 
For the specific operator $B$ corresponding to \eqref{eq:system_linear1} and defined in~\eqref{eq:B}, assumption~\eqref{hyp:B} is satisfied with $q_1=q_2=1$, and an explicit constant $C$ is provided in Appendix~\ref{appendix:A}.

\bigskip





\subsection{Main results and key ideas}


Our work proposes a general approach to address the spectrum of $M$ in \eqref{eq:amplitudes_matrix_short}, leading to a rigorous description of the stability/Turing instability of the original problem \eqref{eq:system}. 
As is common with computed-assisted proofs, we provide sufficient conditions for establishing stability or instability, which can be checked a posteriori, once (most of) the explicit parameters of the system have been fixed. 
The following theorem showcases the kind of results that can be obtained with our approach.

\begin{thm}
\label{thm:main}
    Consider system~\eqref{eq:system} with $\Omega = (0,2)$, $\Omega_1 = (\pi/4, \pi/2)$, $\Omega_2=( \pi/5, \pi/2 + 1/4)$, $a=-3$, $b=2$, $c = 3$, $d=-3$ and $\vartheta =1$. 
    There exists $\delta^*\in [2.428,2.46]$ such that the equilibrium state $(0,0)$ is linearly stable for all $\delta \in [0,\delta^*)$, and linearly unstable for all $\delta \in (\delta^*,4]$ with a single unstable eigenvalue. 
\end{thm}

\begin{remark}
    As will be made clear in the proof, the upper-bound $\delta\leq 4$ in Theorem~\ref{thm:main} is arbitrary, and numerical investigations in fact suggest the system has a single unstable eigenvalue for all $\delta\in (\delta^*,\infty)$. 
    If a statement of the form “the system has a single unstable eigenvalue for all $\delta\in (\delta^{\text{max}},\infty)$" was proven to be true, with an explicit value for $\delta^{\text{max}}$, our approach could in principle be extended to prove Theorem~\ref{thm:main} with $\delta^{\text{max}}$ instead of $4$ as an upper-bound. 
    However, controlling the spectrum (or at least the first eigenvalue) for arbitrarily large values of $\delta$ requires different ideas and techniques than the ones used in this work.
    
    Also, we made no effort to get a particularly sharp enclosure of the transition value $\delta^*$, and tighter bounds could be obtained if needed.
    \end{remark}




\medskip
Let us present here the main steps of our proof of Theorem~\ref{thm:main}, and underline the core ideas. 
Our goal is to control the spectrum of the system matrix $M$ in \eqref{eq:amplitudes_matrix_short} precisely enough to be able to count the number of unstable eigenvalues, i.e., of eigenvalues with positive real part. 

First, we discuss how to control the spectrum for a fixed value of $\delta$. 
Our main tool is the Gershgorin disks theorem. For a \emph{finite} dimensional matrix $M$, it gives a finite union of disks that contains the eigenvalues of the matrix. 
This theorem can be generalized to infinite dimensional operators, under suitable assumptions. 
We provide such a generalization in Section~\ref{sec:Gersh_compact} for the case of infinite dimensional operators having compact resolvent.
However, there are two separate issues that prevent us from successfully applying this theorem directly to $M$: the terms in $M$ coming from the nonlocal operator $B$ are not summable (because for system~\eqref{eq:system_linear1} we have~\eqref{hyp:B} with $q_1=q_2=1$), and there is no reason a priori for the Gershgorin disks to be narrow enough to allow us to conclude regarding the number of unstable eigenvalues. 
We overcome these difficulties by two successive changes of basis. 
The first one takes care of the summability issue, and allows us to prove that eigenvalues of $M$ with large enough index are stable, with an explicit threshold stating what “large enough" means. 
This step holds for any $M$ of the form~\eqref{eq:amplitudes_matrix_short} with $B$ satisfying~\eqref{hyp:B}, and is not computer-assisted. 
The second change of basis is, as it is constructed using finitely many approximate eigenvectors of $M$ computed numerically.
Most of the work then resides in estimating the terms of $M$ in this new basis, in order to get computable upper-bounds for the radii of the Gershgorin disks. 
The bounds we obtain are sharp enough such that, for a fixed $\delta<\delta^*$ not too close to $\delta^*$, we are able to prove that all the Gershgorin disks lie in the half plane $\{Re(z)<0\}$ with negative real part of the complex plane. 
Similarly for a fixed $\delta>\delta^*$ not too close to $\delta^*$, we are able to prove that all but one of the Gershgorin disks lie in the $\{Re(z)<0\}$ part of the complex plane, and that the single remaining Gershgorin disk lies in the $\{Re(z)>0\}$ part of the complex plane. 
Using interval arithmetic~\cite{MooKeaClo09}, we can extend this strategy for all $\delta$ in a small interval, and then repeat it for different intervals of $\delta$ to finally cover the whole interval $[0,4]$.

This whole procedure is successful to prove the instability induced by nonlocality when $\delta$ is large enough, which has been only numerically observed in previous works \cite{reisch2025Turing}. 
However, it does not establish that $\delta^*$ is the  \textit{unique transition point}  from stable to unstable, or \textit{bifurcation point}.
Indeed, when $\delta$ is too close to $\delta^*$, one eigenvalue of $M$ is actually close to zero, and there will eventually be a Gershgorin disk which intersects both $\{Re(z)<0\}$ and $\{Re(z)>0\}$. 
To deal with this critical neighborhood around $\delta^*$ we proceed as follows. We first get a tighter enclosure of the largest eigenvalue $d_0$ using a Newton-Kantorovich type of argument. 
Using essentially the implicit function theorem, we extend these estimates to rigorously control the derivative of $d_0$ with respect to $\delta$, and show that it is positive. This confirms that $d_0(\delta)$ can only cross the imaginary axis once at $\delta^*$, and consequently $\delta^*$ is a bifurcation point for the Turing instability induced by nonlocality.

\medskip

The Gershgorin argument, leading to Theorem~\ref{thm:spectrum}, is presented in details in Section~\ref{sec:Gersh}, whereas the finer analysis around $\delta^*$, culminating in Theorem~\ref{thm:threshold}, is conducted in Section~\ref{sec:threshold}. 
These two theorems taken together immediately imply Theorem~\ref{thm:main}. 
Some of the technical steps of the proofs are presented in the Appendix~\ref{appendix:A}.


\section{Localization of the spectrum.} \label{sec:Gersh}

In this section, we provide a general procedure allowing to describe precisely, for any fixed $\delta$, the spectrum of the \emph{infinite matrix} $M = \widetilde{A} + \delta \widetilde{B}$, with $\widetilde{A}$ the local operator defined in~\eqref{eq:tildeA} and $\widetilde{B}$ the nonlocal operator defined in~\eqref{eq:B1}, assuming $B$ satisfies~\eqref{hyp:B}. In the sequel, we frequently identify (possibly unbounded) linear operators on $L^2(\Omega)\times L^2(\Omega)$ with their representation as infinite matrices in the basis $(\varphi_j)_j$, as was already done in Section~\ref{sec:intro_problem}.

First observe that $M$ has compact resolvent, therefore its spectrum is only composed of eigenvalues.

\begin{proposition} \label{prop:compact_spectrum}
    Let $B$ satisfying \eqref{hyp:B} and let $-\widetilde{A}$ correspond a second-order uniformly elliptic operator with $L^\infty$ coefficients, then $M = \widetilde{A} + \delta \widetilde{B}$ has compact resolvent on $L^2(\Omega)\times L^2(\Omega)$.
\end{proposition}
\begin{proof}
    From the result in \cite[Theorem 3, 6.2]{Evans22} adapted with the Neumann conditions, \cite[V.3]{Luc04}, we have that for any $\lambda \in \C$, with $-\Re(\lambda)$ large enough that $\widetilde{A}-\lambda I : H^2(\Omega)\times H^2(\Omega) \to L^2(\Omega)\times L^2(\Omega)$ is invertible. 
    Since $H^2(\Omega) \hookrightarrow L^2(\Omega)$ is compact then $(\widetilde{A}-\lambda I)^{-1}$ is compact from $L^2(\Omega)\times L^2(\Omega)$ to $L^2(\Omega)\times L^2(\Omega)$. 

    Let $u\in H^2(\Omega)$, it implies $(j^2|u_j|)_{j \in \N}$ is square-summable. Furthermore, from $\eqref{hyp:B}$ with $q_1 > \frac{1}{2}$ and $q_2 > -\frac{3}{2}$, we have by the Cauchy-Schwarz inequality,
    $$ \sum_{i=1}^{+\infty}\left|\sum_{j=1}^{+\infty} B_{i,j}u_j\right|^2 \leq C \sum_{i=1}^{+\infty}i^{-2q_1} \sum_{j=1}^{+\infty}j^{-2(q_2+2)} \sum_{j=1}^{+\infty}(j^2|u_j|)^2 < +\infty,$$ 
    since $-2q_1 < -1$ and $-2(q_2+2) < -1$. That is, $\eqref{hyp:B}$ implies that $\widetilde{B}$ corresponds to bounded operator from $H^2(\Omega)\times H^2(\Omega)$ to $L^2(\Omega)\times L^2(\Omega)$.
    
    Finally, from \cite[Theorem 1.16, IV.1.4]{Kato76}, $M - \lambda I = \widetilde{A} + \delta\widetilde{B} - \lambda I$ is invertible and its inverse is compact from $L^2(\Omega)\times L^2(\Omega)$ to $L^2(\Omega)\times L^2(\Omega)$, since $\delta\widetilde{B}$ is $(-\widehat{A}+\lambda I)$-bounded, for $\lambda \in \C$ with $-\Re(\lambda)$ large enough.
\end{proof}
As a consequence, the spectrum of $M$, $\sigma(M)$, consists of discrete isolated eigenvalues, which we denote $(d_k)_{k\in\N}$ and order by decreasing real part: 
\begin{equation}
    \Re(d_0) \geq \Re(d_1) \geq \dots \geq \Re(d_k) \geq \Re(d_{k+1}) \geq \dots \ . \label{def:eigM}
\end{equation}
In this section, we show how to rigorously enclose these eigenvalues, leading to the following statement in the case of system~\eqref{eq:system}, whose proof is given in Section~\ref{sec:proofspectrum}.
\begin{thm} \label{thm:spectrum}
Consider system~\eqref{eq:system} with $\Omega = (0,2)$, $\Omega_1 = (\pi/4, \pi/2)$, $\Omega_2=( \pi/5, \pi/2 + 1/4)$,
$a=-3$, $b=2$, $c = 3$, $d=-3$ and $\vartheta =1$, and let $\delta_0 = \dfrac{4\times 607}{1000}$ and $\delta_1 = \dfrac{4\times 615}{1000}$. 
We have, 
    \begin{itemize}
        \item[1.] $\forall \delta \in[0,4], \ \exists \mu < 0, \Re(d_1) \leq \mu < \Re(d_0)$, meaning that there is at most one unstable eigenvalue.
        \item[2.] $\forall \delta \in[0,\delta_0]$ $\Re(d_0) < 0$, meaning that system~\eqref{eq:system_linear1} is stable.
        \item[3.] $\forall \delta \in[\delta_1,4]$ $\Re(d_0) > 0$, meaning that system~\eqref{eq:system_linear1} is unstable, with a single unstable eigenvalue.
    \end{itemize}
\end{thm}
Then we can affirm that,
\begin{corollary} \label{cor:d0cont}
    The function $\delta \mapsto d_0(\delta)$ is continuous from $[0,4]$ to $\R$.
\end{corollary}
\begin{proof}
    For any $\delta$, $d_0(\delta)$ is isolated (see \eqref{def:eigM}), and simple (see Theorem~\ref{thm:spectrum}).
    We then have, $d_0(\delta)\in \R$ since $M$ is real.
    Finally, by \cite[IV.3.5]{Kato76}, about the continuity of a system of eigenvalues, we have indeed the continuity of one simple and isolated eigenvalue.
\end{proof}

\begin{figure}[!htbp]
    \centering 
    \centering
        \begin{overpic}[width=0.7\textwidth, trim = 0 100 0 80, clip, grid = false]{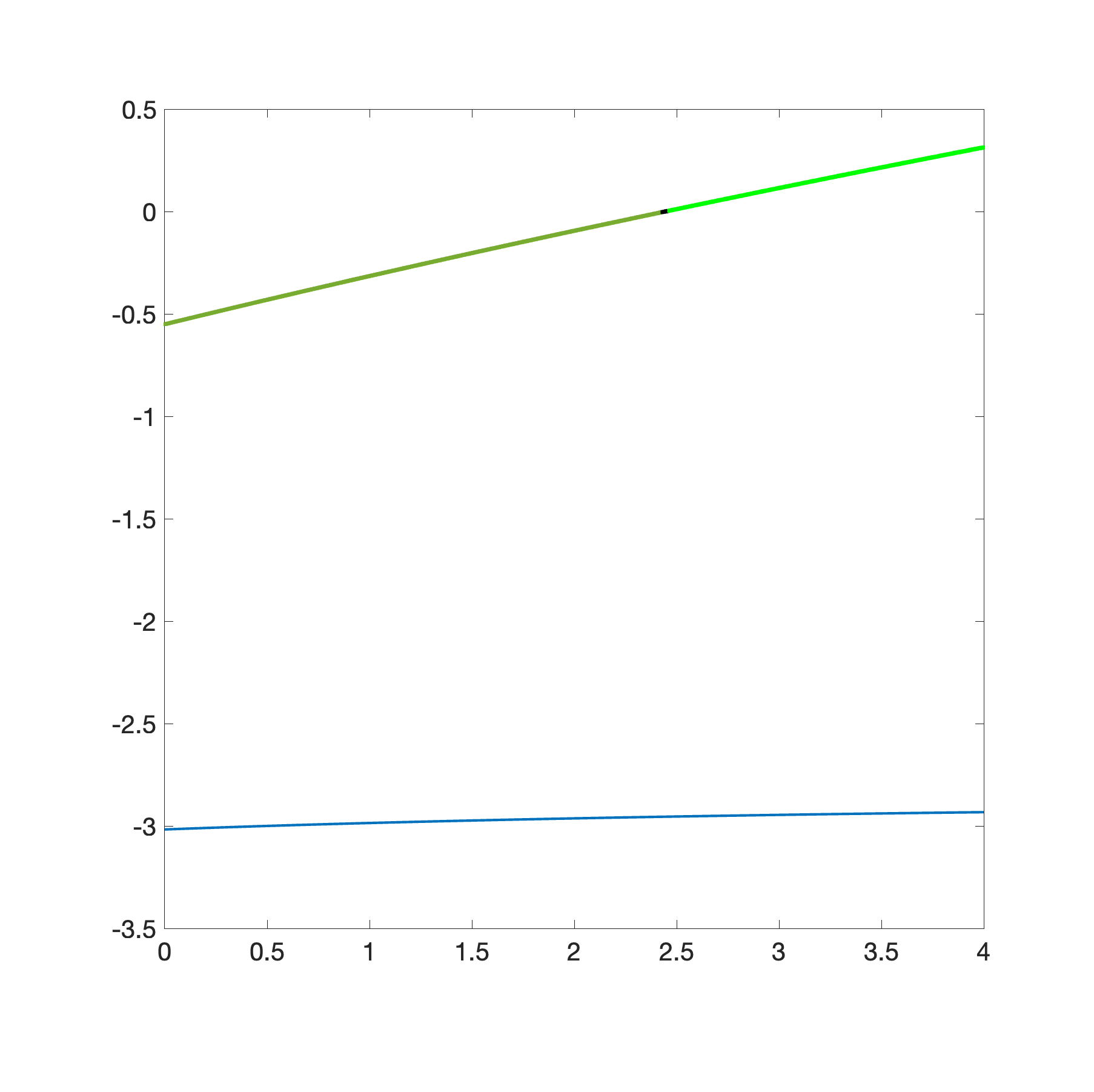}
        \put(550,0){$\delta$}
        \put(120,350){$y$}
        \put(780,160){\color{myblue}\footnotesize $y=\mu(\delta)$}
        \put(400,600){\footnotesize $y = [d_0^-(\delta),d_0^+(\delta)]$}
        \end{overpic}
    \caption{Illustration of some of the quantitative statements obtained within the proof of Theorem~\ref{thm:spectrum} (see Section~\ref{sec:proofspectrum}). 
    We display here a threshold $\mu = \mu(\delta)$ satisfying point 1. of Theorem~\ref{thm:spectrum}, together with a narrow enclosure of $d_0$ given by lower and upper-bounds $d_0^-$ and $d_0^+$, satisfying $d_0^- \leq d_0 \leq d_0^+$ . 
    {\color{darkgreen}{\textbf{Dark green}}} indicates values of $\delta$ for which point $2.$ of Theorem~\ref{thm:spectrum} holds. 
    {\color{lightgreen}{\textbf{Light green}}} indicates values of $\delta$ for which point $3.$ of Theorem~\ref{thm:spectrum} holds. 
    \textbf{Black} indicates values of $\delta$ for which the sign of $d_0$ remains undetermined in Theorem~\ref{thm:spectrum}.}
    \label{fig:d0}
\end{figure}
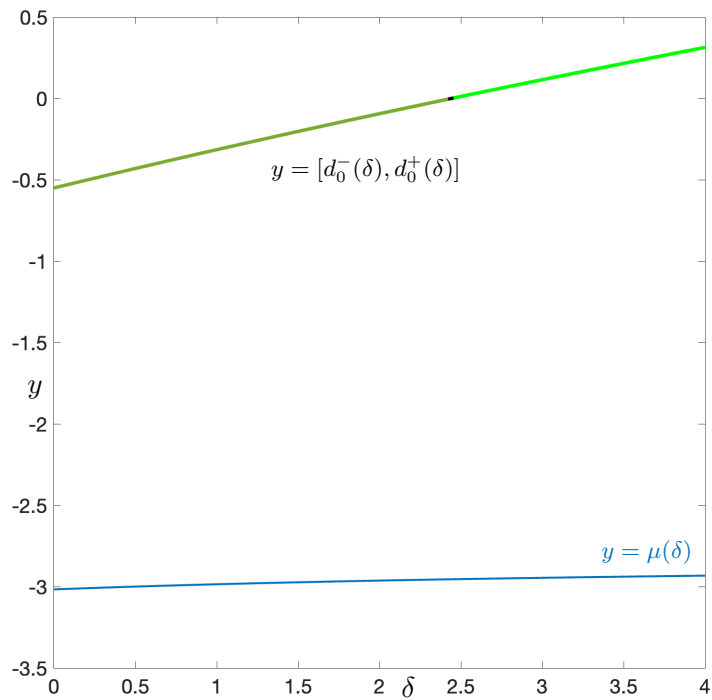

\begin{remark}
Note that $\delta_0=2.428$ and $\delta_1=2.46$, therefore Theorem~\ref{thm:spectrum} implies the linear stability or instability announced in Theorem~\ref{thm:main} when $\delta$ is outside of the interval $(\delta_0,\delta_1)$. 
Within $(\delta_0,\delta_1)$ Theorem~\ref{thm:spectrum} only tells us that there is \emph{at most} one unstable eigenvalue, this case will be further investigated with finer tools in Section~\ref{sec:threshold}.

As will become clear in the remainder of Section~\ref{sec:Gersh}, we prove Theorem~\ref{thm:spectrum} by obtaining quantitative information on the eigenvalues $d_k$ of $M$ that is even more precise than what is stated here. 
For instance, a rigorous enclosure of $d_0$ for all $\delta\in[0,4]$ is provided in Figure~\ref{fig:d0}, together with an explicit value for the threshold $\mu$. 
In fact, the proof of Theorem~\ref{thm:spectrum} relies on a somewhat crude asymptotic estimate for large enough eigenvalues, combined with a much finer control on the first eigenvalues, as illustrated on Figure~\ref{fig:circ1} for $\delta = 1<\delta_0$ and on Figure~\ref{fig:circ2} for $\delta = 4>\delta_1$. 
\end{remark}

\subsection{A Gershgorin theorem for infinite matrices having compact resolvent}
\label{sec:Gersh_compact}

One of the ingredients of the proof of Theorem~\ref{thm:spectrum} is the Gershgorin disk theorem. 
This is a classical result for finite matrices, that can be adapted to some classes of infinite matrices~\cite{HanNetRei68,ShiWilRud87,FarLan91}. 
In order to make the discussion precise, let us introduce some notations and assumptions.
\begin{definition}
     Let $E$ a Banach space having a Schauder basis, and $L:E\to E$ a (possibly unbounded) linear operator, written  $L = (l_{ij})_{(i,j)\in\N^2}$ with respect to the Schauder basis.
     For all $i\in\N$, we denote
     \begin{align*}
        r_i(L) &= \displaystyle \sum_{j\in\N\backslash\{i\}} |l_{ij}|, \\
        \mathcal{D}_i(L) &= \mathcal{D}(l_{ii},r_i(L)) =  \left\{ z\in\C, \vert z-l_{ii}\vert \leq r_i(L) \right\}.
     \end{align*}
     We refer to $\mathcal{D}_i(L)$ as the $i^{th}$ Gershgorin disk of $L$. We note that its radius $r_i(L)$ can be infinite. 
     We also denote by $\sigma(L)$ the spectrum of $L$. 
     In the sequel, we always assume that the space $E$ has the following property
     \begin{align}
     \label{eq:condE}
     \text{for all }x = (x_0, x_1, \dots)\in E, \text{ there exists }i_0\in\N \text{ such that }\sup_{i\in\N}|x_i| = |x_{i_0}|.
         \end{align}
\end{definition}
Under assumption~\eqref{eq:condE}, the classical Gershgorin theorem for finite matrices can easily be generalized to infinite ones, and states that any eigenvalue of $L$ must lie in one of the Gershgorin disks. 
Moreover, for finite matrices there is a stronger version of Gershgorin's theorem, which allows to count eigenvalues within a disjoint subset of Gershgorin disks. 
This stronger version has also been generalized to some infinite matrices in~\cite[Theorem~$2.1$]{FarLan91}. 
However, some of the assumptions of~\cite[Theorem~$2.1$]{FarLan91} are needlessly restrictive (for instance, all the diagonal elements of $L$ have to be nonzero), and others may not be straightforward to check in practice (like the invertibility of a one-parameter family of operators constructed from $L$). 
We propose below a simpler and slightly more general statement, which is still strongly inspired from~\cite[Theorem~$2.1$]{FarLan91} and from its proof, but is easier to use in practice. 
Indeed, in addition to~\eqref{eq:condE}, we simply require that $L$ has compact resolvent. 
Note that, while this assumption is not explicitly made in~\cite[Theorem~$2.1$]{FarLan91}, it is in fact a consequence of the other assumptions of that theorem. 
\begin{thm} \label{thm:Gersh}
    Let $E$ a Banach space having a Schauder basis and satisfying assumption~\eqref{eq:condE}, and $L = (l_{ij})_{(i,j)\in\N^2}: E\to E$ a (possibly unbounded) linear operator.
    Assume that $L$ has a compact resolvent.
    
    Then, the spectrum of $L$, denoted $\sigma(L)$, is included in $\bigcup_{i\in\N} \mathcal{D}_i(L)$.
    Furthermore, if there exist $K\in\N$ and $\{i_1, \dots, i_K\}\subset\N$ such that $\mathcal{R} =\bigcup_{k=1}^{K}\mathcal{D}_{i_k}(L)$ is disjoint from $\bigcup_{i\in\N\backslash\{i_1,\dots,i_K\}}\mathcal{D}_i(L)$, then $\mathcal{R}$ contains exactly $K$ eigenvalues, counted with algebraic multiplicity.
\end{thm}

\begin{proof}
    \textit{Step 1}. Since $L$ has compact resolvent, its spectrum $\sigma(L)$ is formed by a sequence of eigenvalues. For any $\lambda\in\sigma(L)$, consider an associated eigenvector $x$, and $i_0\in\N$ satisfying~\eqref{eq:condE} for that $x$. 
    From $Lx = \lambda x$ we immediately get
    \begin{align*}
        \left(l_{i_0i_0} - \lambda\right) x_{i_0} = \sum_{j\neq i_0} l_{i_0j}x_j,
    \end{align*}
    hence
    \begin{align*}
        \left\vert l_{i_0i_0} - \lambda\right\vert \leq \sum_{j\neq i_0} \left\vert l_{i_0j}\right\vert,
    \end{align*}
    and therefore $\lambda\in \mathcal{D}_{i_0}(L)$.    
    
    \textit{Step 2}. Since $\mathcal{R}$ is disjoint from $\bigcup_{i\in\N\backslash\{i_1,\dots,i_K\}}\mathcal{D}_i(L)$, we can find a closed Jordan curve $\Gamma$ separating $\mathcal{R}$ and the rest of the spectrum of $L$. 
    That is, denoting $U = \text{int}\,\Gamma$, we can find $\Gamma$ such that
    \begin{align*}
        \mathcal{R} \subset U \quad \text{and} \quad \bar{U}\cap \mathcal{D}_i(L) = \emptyset \text{ for }i \in \N\backslash\{i_1,\dots,i_K\}.
    \end{align*}
    Next, we consider $D = \diag(L)$ and the homotopy $L(s) = D + s(L-D)$, for $s \in [0,1]$.
    We have that $\sigma(L(s)) \subset \bigcup_{i\in\N} \mathcal{D}(l_{ii},s r_i(L)) \subset \bigcup_{i\in\N} \mathcal{D}_i(L)$ for all $s \in [0,1]$. 
    Therefore, for each $s$, $\Gamma$ separates $\mathcal{R}$ and the remainder $\bigcup_{i\in\N\backslash\{i_1,\dots,i_K\}}\mathcal{D}_i(L(s))$ of the spectrum of $L(s)$. 
    For $s=0$, $\mathcal{R}$ clearly contains exactly $K$ eigenvalues of the diagonal matrix $L(0)$, and as in the proof of \cite[Theorem~2.1]{FarLan91}, we can then use the semicontinuity of the spectrum of $L$, \cite[Theorem 3.16 IV.3]{Kato76}, to conclude that $\mathcal{R}$ contains exactly $K$ eigenvalues of $L(1) = L$.
\end{proof}
Theorem~\ref{thm:Gersh} shows that Gershgorin disks can be used to control the spectrum of any infinite dimensional operator having compact resolvent (and defined on a space satisfying assumption~\eqref{eq:condE}), which is the case for $M$. 

\subsection{Infinite matrices and basis changes}

We now take a precise look at the operator $M$. 
It can be seen as having $2\times 2$ blocks of elements $M_{2i+\epsilon,2j+\eta}$ for $(i,j) \in \N^2,\, (\epsilon, \eta) \in \{0,1\}^2$  (we start the indexing at $0$ to stay consistent with the numbering of Fourier modes), where
\begin{align*}
M_{2i+\epsilon,2i+\eta} &= \begin{pmatrix}
-\vartheta \lambda_i + a & b \\
c + \delta B_{i,i} & -\lambda_i + d
\end{pmatrix}_{\epsilon,\eta}, \\
M_{2i+\epsilon,2j+\eta} &= \begin{pmatrix}
0 & 0 \\
\delta B_{i,j}  & 0  
\end{pmatrix}_{\epsilon,\eta} , \quad \text{for }i\neq j,
\end{align*}
where $(B_{i,j})$ satisfies \eqref{hyp:B}. 

Naively, if we apply Theorem~\ref{thm:Gersh} to $M$ from \eqref{eq:amplitudes_matrix_short}, we get no information on the localization of the eigenvalues, as some of the disks have infinite radius.
Indeed, for $i\in\N$, the radii $(r_{2i+1}(M))_{i\in\N}$ are not finite since $(B_{i,j})_{j\in\N}$ are not summable.
The purpose of this section is to introduce two successive changes of basis, where the first one ensures that all the radii become finite, and the second one brings enough control on the first disks.

\begin{definition} 
\label{def:spaces}
Let us define \begin{align}
    &\ell^1(\C^2) = \left\{ X=((u_k,v_k))_k \in (\C^2)^\N \ \left| \ \sum_{k=0}^{+\infty} (|u_k|+|v_k|) < +\infty  \right. \right\}, \\
    &\forall\, X\in\ell^1(\C^2), \ \|X\|_1 = \sum_{k=0}^{+\infty} (|u_k|+|v_k|).
\end{align}
$(\ell^1(\C^2), \|\cdot\|_1)$ is a Banach space, the canonical basis of $\C^\N$ is a Schauder basis, and if $X\in\ell^1(\C^2), \ \sup\left( \, \{|u_k|,\, k\in\N\} \cup \{|v_k|,\, k\in\N\}\right)$ exists and is reached, therefore $(\ell^1(\C^2), \|\cdot\|_1)$ satisfies assumption~\eqref{eq:condE}.
\end{definition}

By elliptic regularity, we know that eigenvectors of $M$ belong to $\ell^1(\C^2)$. From now on and until the end of Section~\ref{sec:Gersh}, infinite matrices correspond to (bounded or unbounded) linear operators on $\ell^1(\C^2)$. Note that the specific choice of sequence space is irrelevant in this section, provided assumption~\eqref{eq:condE} holds and all the eigenvectors of $M$ belong to that space, therefore one could also have considered some $\ell^2$ space here.

We now introduce the first change of basis, which will allow us to recover finite Gershgorin disks.
\begin{definition}
\label{def:Q}
Let $f : x \mapsto \max(1,x^p)$, $p\geq 0$, and $Q$ be the infinite diagonal matrix such that
\begin{align*}
    Q_{2i+\epsilon,2j+\eta} = \frac{1}{f(i)}\delta_{2i+\epsilon,2j+\eta}, \quad\text{where }\delta_{k,l} = \begin{cases}
1, \text{ if } k=l \\
0 \text{ otherwise}
\end{cases}.
\end{align*} 

We have $Q \in \mathcal{L}(\ell^1(\C^2))$, and $Q$ is bounded. 
Furthermore, $Q$ is invertible but its inverse is unbounded when $p>0$. 

We then define $\widetilde{M} = Q^{-1}MQ$. 
\end{definition}

\begin{remark}
    In principle, one could use different choices of $f$. 
    The above choice proved sufficient for establishing Theorem~\ref{thm:main}, but for more involved nonlocal operators $B$ a different $f$ might prove more efficient.

    In order to obtain $\widetilde{M}$ from $M$, one simply has to multiply the rows of index $2i^{th}$ and $2i+1^{th}$ by $f(i)$, and the columns of index $2j^{th}$ and $2j+1^{th}$ by $\frac{1}{f(j)}$, for all $i,j\geq 0$:
    \begin{equation*}
    \widetilde{M} = \begin{pmatrix}
        a-\vartheta\lambda_0 & b &  0 & 0 & 0 & 0 & \dots \\
        c + \delta B_{00} & d- \lambda_0 & \delta \frac{f(0)}{f(1)} B_{01} & 0 & \delta \frac{f(0)}{f(2)} B_{02} & 0 &\dots \\
        0 & 0 & a-\vartheta\lambda_1 & b & 0 & 0 &\ddots \\
        \delta \frac{f(1)}{f(0)} B_{10} & 0 & c + \delta B_{11} & d-\lambda_1 & \delta \frac{f(1)}{f(2)} B_{12} & 0 & \ddots\\
        0 & 0 & 0 & 0 & a-\vartheta\lambda_2 & b & \ddots \\
        \delta \frac{f(2)}{f(0)} B_{20} & 0 & \delta \frac{f(2)}{f(1)} B_{21} & 0 & c + \delta B_{22} & d-\lambda_2 & \ddots \\
        \vdots & \vdots  & \ddots  & \ddots & \ddots  & \ddots & \ddots \\
    \end{pmatrix}.
    \end{equation*}
\end{remark}

By choosing an appropriate $p$ in Definition~\ref{def:Q}, $\widetilde{M}$ will have Gershgorin disks of finite radius. Furthermore, because of the asymptotic behavior of the eigenvalues $(\lambda_i)_{i\in\N}$ of the negative Laplacian, these disks $\mathcal{D}_i(\widetilde{M})$ will all be contained in $\{\Re(z)< 0\}$ for $i$ large enough.
Indeed, by Weyl's law, we know that there exists a constant $\varkappa>0$ such that
\begin{equation}\label{Weyl}
    \lambda_i \geq \varkappa i^{\frac{2}{n}}, \quad \forall i \in \N.
\end{equation}
When $p$ is in an appropriate range (depending on $q_1,q_2$ from~\eqref{hyp:B} and on the dimension $n$), we can then give an explicit threshold $i_0$ (depending on the constant $\varkappa$ from~\eqref{Weyl}) after which all the Gershgorin disks of $\widetilde{M}$ are in the left half of the complex plane.


\begin{lem}\label{lem:new}
    Let $p\in(1-q_2,q_1+\frac{2}{n})$, and $i_0\in\N$ such that, for all $i\geq i_0$,
    \begin{align*}
        -\vartheta \varkappa i^{\frac{2}{n}} + a + |b| < 0 \quad \text{and} \quad 
       -  \varkappa i^{\frac{2}{n}} + C'\vert\delta\vert i^{p-q_1} + d + |b| < 0, 
    \end{align*}
    with $\varkappa$ from~\eqref{Weyl} and $C' = C\sum_{j\geq 0}1/\max(1,j^{p+q_2})$ with $C$ from \eqref{hyp:B}.
    Then
    \begin{equation}\label{inclusion}
    \bigcup_{k\geq 2i_0}\mathcal{D}_k(\widetilde{M})\subset \left\{ z\in\C, \ \Re(z) < 0 \right\}.
    \end{equation} 
\end{lem}
\begin{remark}
    For many specific cases of the domain $\Omega$, say a cube $\Omega = (0,l)^n$ or a sphere $\Omega = \{|x| < l\}$, the constant $\varkappa$ in \eqref{Weyl} can be explicitly estimated. In those cases, one can easily derive from Lemma \ref{lem:new} an explicit value of $i_0$ for which, the $k$-th Gershgorin disk lies completely on the left half plane for all $k\geq 2i_0$. See also Lemma~\ref{lem:last_disks} for explicit bounds, in a slightly more complicated setting.
\end{remark}
\begin{proof}
    For all $i\in \N$,
    $r_{2i}(\widetilde{M}) = |b|$ and 
    \begin{equation}\label{eq:roddMtilde}
        r_{2i+1}(\widetilde{M}) \leq |\delta|\sum_{j=0}^{+\infty} \frac{f(i)}{f(j)}|B_{i,j}|+  |c|.
    \end{equation}
    Using \eqref{hyp:B} and $f(k) = \max(1,k^p) $, we get
    \begin{align}
        r_{2i+1}(\widetilde{M}) &\leq C|\delta| i^{p-q_1}\sum_{j=0}^{+\infty} \frac{1}{\max(1,j^{p+q_2})}| +  |c| \nonumber\\
        &\leq C'|\delta| i^{p-q_1} + |c|, \label{eq:revenMtilde}
    \end{align}
    since $p+q_2 > 1$. Moreover, for all $i\in\N$,
    \begin{align*}
        \widetilde{M}_{2i,2i} = M_{2i,2i} = -\vartheta \lambda_i + a \qquad \text{and} \qquad \widetilde{M}_{2i+1,2i+1} = M_{2i+1,2i+1} = - \lambda_i + d.
    \end{align*}
    Therefore,
    \begin{align}\label{eq:bound_even}
        \widetilde{M}_{2i,2i} + r_{2i}(\widetilde{M}) \leq -\vartheta \varkappa i^{\frac{2}{n}} + a + |b|,
    \end{align}
    and
    \begin{align}\label{eq:bound_odd}
        \widetilde{M}_{2i+1,2i+1} + r_{2i+1}(\widetilde{M}) \leq - \varkappa i^{\frac{2}{n}} + C'\vert\delta\vert i^{p-q_1} + d + |b|.
    \end{align}
    Since $\vartheta,\varkappa,C'>0$ and $p-q_1<2/n$, the right hand side of~\eqref{eq:bound_even} and of~\eqref{eq:bound_odd} is negative for all $i$ large enough, hence there exists $i_0$ satifying the assumptions of the lemma, and for all $k\geq 2i_0$, $\mathcal{D}_k(\widetilde{M})\subset \left\{ z\in\C, \ \Re(z) < 0 \right\}$.
\end{proof}

Lemma \ref{lem:new} already shows for $\widetilde{M}$ that the $k$-th Gershgorin disks lie on the left half plane for large $k$. 
However, for small $k$ we do not yet have enough control to be able to determine the sign of the real parts of the eigenvalues of $\widetilde{M}$. 
Therefore, we are going to also approximately diagonalize a finite submatrix of $\widetilde{M}$. 
To that end, we first introduce the truncation operator $\Pi^N$, for $N\in\N$.

\begin{definition}
    Let $N\in\N$, we denote with $\Pi^N$ the following infinite matrix,
    \begin{equation}
        \forall (i,j)^2 \in \N^2, \, (\epsilon,\eta) \in \{0,1\}^2, \ (\Pi^N)_{2i+\epsilon,2j+\eta} = \begin{cases}
            \delta_{2i+\epsilon,2j+\eta}, \ i,j < N \\
            0, \ \text{otherwise}.
        \end{cases}
    \end{equation}
    \begin{equation*}
        \Pi^N = \begin{pmatrix}
        1 & 0 & 0 & \cdots & \cdots & 0 & \cdots \\
        0 & 1 & 0 & \ddots &  & 0 & \cdots \\
        0 & 0 & 1 & \ddots & \ddots & 0 & \\
        \vdots  & \ddots  & \ddots  & \ddots & \ddots  &  \vdots &  \\
        \vdots  &  & \ddots & \ddots & 1 & 0 &\ddots \\
        0 & 0 & 0 & \cdots & 0 & 0 &\ddots \\
        \vdots & \vdots  &  &  & \ddots  & \ddots & \ddots
    \end{pmatrix} 
    \begin{array}{ll}  
    \xleftarrow{}0^{th} \ \text{row} \phantom{\cdots}\\
    \phantom{\ddots}\\ 
    \phantom{\ddots}\\ 
    \phantom{\vdots}\\ 
    \phantom{\vdots}\\ 
    \xleftarrow{}2N^{th} \ \text{row} \phantom{\vdots}\\ 
    \phantom{\vdots}
    \end{array}
    \end{equation*}
$\Pi^N$ is the canonical projection from $(\C^2)^{\N}$ into $(\C^2)^N$. 
\end{definition}

\begin{remark}
In the sequel, given an infinite matrix $L$, we frequently identify $\Pi^N L \Pi^N$ with a finite $2N\times 2N$ matrix. 
Reciprocally, we often identify a given $2N\times 2N$ matrix $L_N$ with an infinite matrix obtained by completing $L_N$ with zeros.
\end{remark}

\begin{definition} \label{def:P}
    Let $P_N$ be an invertible matrix of size $2N\times 2N$, obtained numerically such that $P_N^{-1}\, \Pi^N \widetilde{M}\Pi^N\, P_N$ is approximately diagonal (that is, the columns of $P_N$ are taken to be numerically computed approximate eigenvectors of $\Pi^N \widetilde{M}\Pi^N$).  
    We then define $P$, an infinite matrix,  by
    \begin{equation}
        P = P_N + I - \Pi^N,
    \end{equation}
    which can schematically be represented as follows:
    \begin{equation*}
        P = 
    \begin{pmatrix}
    \begin{array}{c|} P_N \\ \hline \end{array} & \hphantom{P_N}  \\
    \hphantom{P_N} & I 
    \end{pmatrix}.
    \end{equation*}
    Note that $P$ is invertible in $\mathcal{L}(\ell^1(\C^2))$, and its inverse is bounded, $P^{-1} = (P_N)^{-1} + (I-\Pi^N)$.

    Finally, we consider 
\begin{equation}
\label{eq:deffrakM}
    \M = P^{-1}\widetilde{M}P,
\end{equation}
which is an unbounded operator on $\ell^1(\C^2)$.
\end{definition}
By construction, $\M$ has the same spectrum as $M$, and we are going to apply Theorem~\ref{thm:Gersh} to $\M$.

\begin{remark}
    Thanks to our choice of $P$, $\Pi^N\M\Pi^N$ should be “almost" diagonal, in the sense that we expect to have $\M_{2i+\epsilon,2j+\eta} \approx 0$, for $2i+\epsilon \neq 2j+\eta$ with $i,j < N$, $\epsilon,\eta \in \{ 0, 1\}$.
\end{remark}
 
\subsection{Gershgorin disks and estimation of radius bounds}

Our goal is now to estimate the radii of the Gershgorin disks of $\M$ defined in~\eqref{eq:deffrakM}, in order to precisely localize the spectrum of $M$ and to prove Theorem~\ref{thm:spectrum}.

\begin{proposition} \label{prop:radiibounds}
    Consider $B$ and $q_1,q_2$ as in~\eqref{hyp:B}, $p>1-q_2$, $P$ as in Definition~\ref{def:P} and Q as in Definition~\ref{def:Q}, and $\M$ as in~\eqref{eq:deffrakM}, i.e., $\M = P^{-1}Q^{-1} M Q P$. 
    Then, $\M$ satisfies the hypotheses of Theorem~\ref{thm:Gersh}, and we have the following estimates on the radii of its Gershgorin disks:
    \begin{align} 
    \intertext{When $i<N$ and $\epsilon \in \{0,1\}$}
        r_{2i + \epsilon}(\M)
        &\leq \dfrac{ C|\delta| }{(p+q_2-1)(N-1)^{p+q_2-1}} \sum_{k=0}^{N-1}\dfrac{|P_{2i+\epsilon,2k+1}^{-1}|}{\max(1,k)^{-p+q_1}}  + \sum_{j=0}^{N-1} \underset{ \eta \neq \epsilon \text{ if } j=i}{\sum_{\eta = 0,}^{1}}  |\M_{2i+\epsilon,2j+\eta}| \notag \\
        &:= R_{2i+\epsilon}.  \label{eq:rd1}
    \end{align}
    \begin{align}
    \intertext{When $i\geq N$} 
        r_{2i}(\M) &= |b| := R_{2i}, \label{eq:rd2} \\ 
        r_{2i+1}(\M) &\leq C|\delta|i^{p-q_1}\left(\sum_{j=0}^{N-1} \sum_{k=0}^{N-1} \dfrac{|P_{2k,2j}| + |P_{2k,2j+1}|}{\max(1,k)^{p+q_2}} + \dfrac{1}{(p+q_2-1)(N-1)^{p+q_2-1}}\right) + |c| \notag \\
        &:= R_{2i+1}. \label{eq:rd3} 
    \end{align}
\end{proposition}

\begin{proof}
    We derive upper-bounds for the quantities $r_{2i+\epsilon}(\M)$, for all $i\in \N$ and $\epsilon\in\{0,1\}$ by splitting them into four cases. 
    Let $i\in\N, \epsilon \in\{0,1\}$ and let $j\in \N$, $\eta\in\{0,1\}$.

\begin{enumerate}
    \item When $i<N$ and $j<N$, this is the only case where we do not actually estimate but in fact compute $\M_{2i+\epsilon,2j+\eta}$ explicitly on the computer (rigorously using interval arithmetic). 
    Thanks to the choice of $P$, we “almost" diagonalize $\Pi^N\widetilde{M}\Pi^N$. 
    We expect to have
    \begin{align} \label{eq:coeff1} \M_{2i+\epsilon,2j+\eta} \approx \delta_{2i+\epsilon,2j+\eta} d_{2i+\epsilon}, \end{align}
    where the $(d_k)_{k\in\N}$ are the eigenvalues of $M$ \eqref{def:eigM}.
    \item When $i<N$ and $j\geq N$, gathering the definitions of $P$, $Q$ and $M$, we have
    \begin{align} \label{eq:coeff2}
    \M_{2i+\epsilon,2j+\eta} &= \sum_{k=0}^{2N-1} P^{-1}_{2i+\epsilon,k}\widetilde{M}_{k,2j+\eta}, \notag \\
    &= \dfrac{1}{f(j)}\sum_{k=0}^{2N-1} f\left(\left\lfloor \frac{k}{2} \right\rfloor\right) P^{-1}_{2i+\epsilon,k}M_{k,2j+\eta} \notag \\
    & =(1-\eta)\dfrac{\delta}{f(j)}\sum_{k=0}^{N-1} f(k) P^{-1}_{2i+\epsilon,2k+1}B_{k,j}.
    \end{align}
    \item When $i\geq N$ and $j<N$, we have in the same way,
    \begin{align} \label{eq:coeff3}
    \M_{2i+\epsilon,2j+\eta} &= \epsilon \delta f(i) \sum_{k=0}^{N-1} \frac{1}{f(k)} B_{i,k}P_{2k,2j+\eta}.
    \end{align}
    \item When $i\geq N$ and $j\geq N$,
    \begin{align} \label{eq:coeff4}
    \M_{2i+\epsilon,2j+\eta} &= \widetilde{M}_{2i+\epsilon,2j+\eta} \notag \\
    &= \epsilon(1-\eta)  \delta\dfrac{f(i)}{f(j)} B_{i,j} + \delta_{ij}\begin{pmatrix} -\vartheta\lambda_{ii} + a & b \\
    c & -\vartheta\lambda_{ii} + a
    \end{pmatrix}_{\epsilon, \eta}.
    \end{align}
\end{enumerate}
    We can now deduce bounds on the radii.
    Thanks to \eqref{eq:coeff2}, for $i<N$ and $\epsilon \in\{0,1\}$ we have
    \begin{align*}
        r_{2i+\epsilon}(\M) &= \sum_{j=0}^{+\infty} \underset{ \eta \neq \epsilon \text{ if } j=i}{\sum_{\eta = 0,}^{1}} |\M_{2i+\epsilon,2j+\eta}| \\
        &= \sum_{j=0}^{N-1} \underset{ \eta \neq \epsilon \text{ if } j=i}{\sum_{\eta = 0,}^{1}} |\M_{2i+\epsilon,2j+\eta}|  + \sum_{j=N}^{+\infty} \sum_{\eta = 0}^{1} |\M_{2i+\epsilon,2j+\eta}| \\
        &= \sum_{j=0}^{N-1} \underset{ \eta \neq \epsilon \text{ if } j=i}{\sum_{\eta = 0,}^{1}} |\M_{2i+\epsilon,2j+\eta}|  + \sum_{j=N}^{+\infty}\dfrac{|\delta|}{f(j)} \left| \sum_{k=0}^{N-1}f(k)P^{-1}_{2i+\epsilon,2k+1}B_{k,j} \right| \\
        &\leq \sum_{j=0}^{N-1} \underset{ \eta \neq \epsilon \text{ if } j=i}{\sum_{\eta = 0,}^{1}} |\M_{2i+\epsilon,2j+\eta}|  + C|\delta| \sum_{k=0}^{N-1}\dfrac{f(k)|P^{-1}_{2i+\epsilon,2k+1}|}{\max(1,k^{q_1})} \sum_{j=N}^{+\infty}\dfrac{1}{f(j)j^{q_2}}.
    \end{align*}

Using that $f(j) = j^p$ for $j\geq N$, and Lemma~\ref{lem:inequality} in Appendix~\ref{appendix:A}, we obtain~\eqref{eq:rd1}. 

Next, thanks to \eqref{eq:coeff3}, \eqref{eq:coeff4}, for $i\geq N$ and $\epsilon \in\{0,1\}$ we get
    \begin{align*}
        r_{2i+\epsilon}(\M) &= \sum_{j=0}^{+\infty} \underset{ \eta \neq \epsilon \text{ if } j=i}{\sum_{\eta = 0,}^{1}} |\M_{2i+\epsilon,2j+\eta}| \\
        &= \sum_{j=0}^{N-1}  \sum_{\eta = 0}^{1} |\M_{2i+\epsilon,2j+\eta}|  + \sum_{j=N}^{+\infty} \underset{ \eta \neq \epsilon \text{ if } j=i}{\sum_{\eta = 0,}^{1}} |\M_{2i+\epsilon,2j+\eta}| \\
        &= \epsilon |\delta| f(i) \sum_{j=0}^{N-1}  \sum_{\eta = 0}^{1}\left|\sum_{k=0}^{N-1} \dfrac{1}{f(k)} B_{i,k} P_{2k,2j+\eta} \right| \\
        &\quad + \epsilon \left(|\delta|f(i)\sum_{j=N,j\neq i}^{+\infty}\dfrac{|B_{i,j}|}{f(j)} + |c+\delta B_{i,i}| \right) + (1-\epsilon) |b| \\
        &\leq \epsilon C|\delta|\dfrac{f(i)}{i^{q_1}} \left( \sum_{j=0}^{N-1} \sum_{k=0}^{N-1} \dfrac{|P_{2k,2j}|+ |P_{2k,2j+1}|}{f(k)\max(1,k^{q_2})} + \sum_{j=N}^{+\infty} \dfrac{1}{f(j)j^{q_2}} \right) + \epsilon|c| + (1-\epsilon)|b|. \\
    \end{align*}
Using once again Lemma~\ref{lem:inequality}, we obtain~\eqref{eq:rd2} when $\epsilon=0$, and~\eqref{eq:rd3} when $\epsilon=1$.
\end{proof}
The key observation is that all the quantities $R_{2i+\epsilon}$ for $i<N$ and $\epsilon\in\{0,1\}$ only involve finite computations, and therefore can be obtained explicitly. 
Similarly, the corresponding diagonal elements of $\M$, i.e., the centers of the Gershgorin disks, can be obtained explicitly. 
Therefore, we have an explicit control on the localization of the first $2N$ disks. Note that, in order to establish Theorem~\ref{thm:spectrum}, we do not need a very precise control on $\mathcal{D}_i(\M)$ for $i\geq 2N$, and so we are simply going to get a uniform in $i$ estimate, showing that all these disks lie to the left of some vertical line $\{\Re(z) = \mu\},\, \mu < 0$, in the complex plane.
Provided $N$ is taken large enough, this behavior is to be expected, as shown by Lemma~\ref{lem:new}. However, note that we cannot simply apply this lemma here, as the change of basis $P$ that transformed $\widetilde{M}$ into $\mathfrak{M}$ affects the radius of the Gershgorin disks for $i\geq 2N$. Nonetheless, the argument to come in Lemma~\ref{lem:last_disks} is similar in spirit with Lemma~\ref{lem:new}. 

We now restrict our attention to one-dimensional domains of the form $\Omega = (0,l)$, like in Theorem~\ref{thm:main}, for which the eigenvalues of the negative Laplacian operator are explicitly known: 
$$
\lambda_i = \left( \dfrac{\pi i}{l} \right)^2,\quad  \text{for all }i\in\N,
$$
which will allow us to get fully computable bounds.
We first introduce some notation, and then provide the necessary control on the disks $\mathcal{D}_i(\M)$ for $i\geq 2N$ in Lemma~\ref{lem:last_disks}.

\begin{definition} \label{def:centers}
For all $i\geq N$ and $\epsilon\in\{0,1\}$, we denote by $c_{2i+\epsilon}$ the center of the $2i+\epsilon$-th Gershgorin disk. That is: 
\begin{align*}
    c_{2i} &= \M_{2i,2i} =  - \vartheta \left(\dfrac{\pi i}{l}\right)^2 + a, \\
    c_{2i+1} &= \M_{2i+1,2i+1} =  -\left(\dfrac{\pi i}{l}\right)^2 + d. 
\end{align*}
We also recall the corresponding radii bounds computed in Proposition~\ref{prop:radiibounds}, for all $i\geq N$:
\begin{align}
    R_{2i} &= |b|, \nonumber \\
    R_{2i+1} &= C|\delta|\rho_{N,p} i^{p-q_1} + |c|, \label{eq:roddMfrak}
    \intertext{with,}
    \rho_{N,p} &= \displaystyle \sum_{j=0}^{N-1} \sum_{k=0}^{N-1} \dfrac{|P_{2k,2j}| + |P_{2k,2j+1}|}{\max(1,k)^{p+q_2}} + \dfrac{1}{(p+q_2-1)(N-1)^{p+q_2-1}}. \nonumber
\end{align}
Finally, we define
\begin{align*}
    \overline{m}_{N,p} &= \sup_{k \geq 2N} \left( c_k + R_k \right) .
\end{align*}
\end{definition}


Since $\bigcup_{k \geq 2N} \mathcal{D}_k(\M) \subset \{z \in \C \ |\ \Re(z) \leq \overline{m}_{N,p} \} $, our remaining task is to get a computable upper-bound for $\overline{m}_{N,p}$.

\begin{remark}
For a fixed $N$, and any $i\geq N$, the formula~\eqref{eq:roddMfrak} for $R_{2i+1}$ scales like $i^{p-q_1}$ (formal calculations show that $P_{2k,2j}$ is proportional to $k^p$, hence  $\sum\limits_{j=0}^{N-1} \sum\limits_{k=0}^{N-1} \dfrac{|P_{2k,2j}| + |P_{2k,2j+1}|}{\max(1,k)^{p+q_2}}$ does not depend on $p$). 
Therefore, we need to take $p\leq q_1+2$ to ensure that the radii do no grow faster than the eigenvalues of the Laplacian, and we would in fact like to take $p$ as small as possible, in order to get a better control of the radii for $i\geq N$. 
However, we must at the same time have $p>1-q_2$, otherwise $R_{2i+1}$ becomes infinite. Moreover, for $i< N$, the formula~\eqref{eq:rd1} for $R_{2i+1}$ is proportional to $\frac{1}{N^{p+q_2-1}}$, therefore taking $p$ too close to $1-q_2$ is going to be detrimental for the control of first radii. For our concrete example (for which $q_1=q_2=1$) taking $p$ close to $q_1+1 = 2$ proved to be a suitable compromise, but we do not claim that this is the optimal choice in general, and therefore provide below a computable value of $\overline{m}_{N,p}$ for all $p\in(1-q_2,q_1+2]$.
\end{remark}

\begin{lem} \label{lem:last_disks}
Take $\Omega = (0,l)$ in~\eqref{eq:system_linear1}, and let $\delta \neq 0$, $N \in \N$, $N > 1$. 
Let $p>1-q_2$,
\begin{itemize}
    \item if $p\in(1-q_2,q_1]$, then $\overline{m}_{N,p} = \displaystyle\max_{\epsilon\in\{ 0, 1\}} (c_{2N+\epsilon} + R_{2N+\epsilon})$;
    \item  if $p\in(q_1,q_1+2)$, then $\overline{m}_{N,p} = \max(c_{2N} + R_{2N},c_{2N_1+1} + R_{2N_1+1})$, where 
    $$
    N_1 = \max(N, \underset{i \in \{ \lfloor \bar{a}_{N,p} \rfloor, \lceil  \bar{a}_{N,p} \rceil \}}{\arg\max} (c_{2i+1}+R_{2i+1})),
    $$
    with 
    $$\bar{a}_{N,p} = \left(\dfrac{C|\delta| l^2\rho_{N,p}(p-q_1)}{2\pi^2}\right)^{\frac{1}{2+q_1-p}};
    $$
    \item if $p=q_1+2$ and $C|\delta|\rho_{N,p} \leq \left( \dfrac{\pi}{l}\right)^2$, then $\overline{m}_{N,p} = \displaystyle\max_{\epsilon\in\{ 0, 1\}} (c_{2N+\epsilon} + R_{2N+\epsilon})$;
    \item otherwise, $\overline{m}_{N,p} = + \infty$.
\end{itemize}
\end{lem}

\begin{remark}
    For $\delta = 0$, $\overline{m}_{N,p} = \displaystyle\max_{\epsilon\in\{ 0, 1\}} (c_{2N+\epsilon} + R_{2N+\epsilon})$, $\overline{m}_{N,p}$ no longer depends on $B$ or $p$.
\end{remark}

\begin{remark}
    Compared to Lemma~\ref{lem:new}, one difference with the situation considered in Lemma~\ref{lem:last_disks} is that we now have to deal with the extra change of basis $\M = P^{-1} \widetilde{M} P$, which leads to the additional term $\rho_{N,p}$ (compare~\eqref{eq:roddMtilde} and~\eqref{eq:roddMfrak}) depending on the truncation parameter $N$. 
    Moreover, whereas in Lemma~\ref{lem:new} we only proved the \emph{existence} of a threshold after which all Gershgorin disks were in the left half-plane, here we need an \emph{explicit} control on all the disks $\mathcal{D}_i(\M)$ for $i\geq 2N$ in order to obtain a value of $\mu$ in Theorem~\ref{thm:spectrum}. 
    Which is why we have to work slightly more in the proof of Lemma~\ref{lem:last_disks} than in that of Lemma~\ref{lem:new}.

\end{remark}

\begin{proof}
    Let $p > 1-q_2 $. For $\epsilon = 0$, since $\vartheta \geq 0$, $\displaystyle \sup_{i\geq N} (c_{2i} + R_{2i}) = c_{2N} + R_{2N}  = - \vartheta \left( \dfrac{N \pi}{l} \right)^2 + a + |b|$. For $\epsilon = 1$ and $i\geq N$, we have
\begin{align*}
    c_{2i+1} + R_{2i+1} = - \left(\dfrac{\pi}{l}\right)^2 i^2 + C|\delta|\rho_{N,p} i^{p-q_1} + |c| + d .
\end{align*}
Now we consider the following cases
\begin{itemize}
\item 
If $p > q_1+2$, then $c_{2i+1} + R_{2i+1}$ increases with $i$ for $i$ large enough, therefore $\displaystyle \sup_{i \geq N} c_{2i+1} + R_{2i+1} = +\infty$.

\item If $p = q_1+2$ and $C|\delta|\rho_{N,p} > \left(\dfrac{\pi}{l} \right)^2$, then $c_{2i+1} + R_{2i+1}$ increases with $i$, therefore $\displaystyle \sup_{i \geq N} c_{2i+1} + R_{2i+1} = +\infty$.

\item If  $p = q_1+2$ and $C|\delta|\rho_{N,p} \leq \left(\dfrac{\pi}{l} \right)^2$, then $c_{2i+1} + R_{2i+1}$ decreases with $i$, therefore $\displaystyle \sup_{i \geq N} c_{2i+1} + R_{2i+1} =  c_{2N+1} + R_{2N+1}$.

\item We are left with the case $p < q_1+2$. 
Let us denote $g : x \mapsto -\left(\dfrac{\pi}{l}\right)^2 x^2 + C|\delta|\rho_{N,p} x^{p-q_1} + |c| + d$. The function $g$ is differentiable on $\R_+^*$ and 
$$\forall x \in \R_+^*, \ g'(x) = \left[ C|\delta|\rho_{N,p}(p-q_1) x^{p-q_1-2} - 2\left(\dfrac{\pi}{l}\right)^2\right]x.$$
There are two subcases
\begin{itemize}
    \item[$\circ$]
If $p \leq q_1$, $\forall x \in\R_+^*, \ g'(x) \leq 0$. Then, $g$ is decreasing, therefore $\displaystyle \sup_{i \geq N} c_{2i+1} + R_{2i+1} =  c_{2N+1} + R_{2N+1}$.

\item[$\circ$] If $p > q_1$, with $\bar{a}_{N,p} = \left(\dfrac{C|\delta| l^2\rho_{N,p}(p-q_1)}{2\pi^2}\right)^{\frac{1}{2+q_1-p}}$,  $g'(\bar{a}_{N,p}) = 0$ and $g(\bar{a}_{N,p}) = \sup_{x\in\R^*_+} g(x)$. Therefore, $\displaystyle \sup_{i\geq N} (c_{2i+1} + R_{2i+1}) = c_{2N_1+1} + R_{2N_1+1}$ where \\ $N_1 = \max(N,  \underset{i \in \{ \lfloor \bar{a}_{N,p} \rfloor, \lceil  \bar{a}_{N,p} \rceil \}}{\arg\max} (c_{2i+1}+R_{2i+1}))$. 
\qedhere
\end{itemize}
\end{itemize}

\end{proof}

\subsection{Proof of Theorem~\ref{thm:spectrum}}
\label{sec:proofspectrum}

In the previous subsection, we have obtained a precise control on the $2N$ first Gershgorin disks of $\M$, as well as some rougher but still explicit estimates on the others disks. 
We now use these results in order to prove Theorem~\ref{thm:spectrum} on the location of the spectrum of $M$.

The proof uses interval arithmetic with the library \textsc{Intlab} from \cite{intlabRump}. 
Interval arithmetic not only allows us to rigorously control rounding errors, but it also enables us to derive estimates on the location on the Gershgorin disks which are valid for all $\delta$ in a (relatively small) interval.
The proof is presented just below and its computational parts can be reproduced using the code available at~\cite{GitHub}.

\begin{figure}[!htbp]
    \centering
    \begin{subfigure}{\linewidth}
    \centering
    \begin{overpic}[width=0.75\linewidth, trim = 0 50 0 30, grid = false]{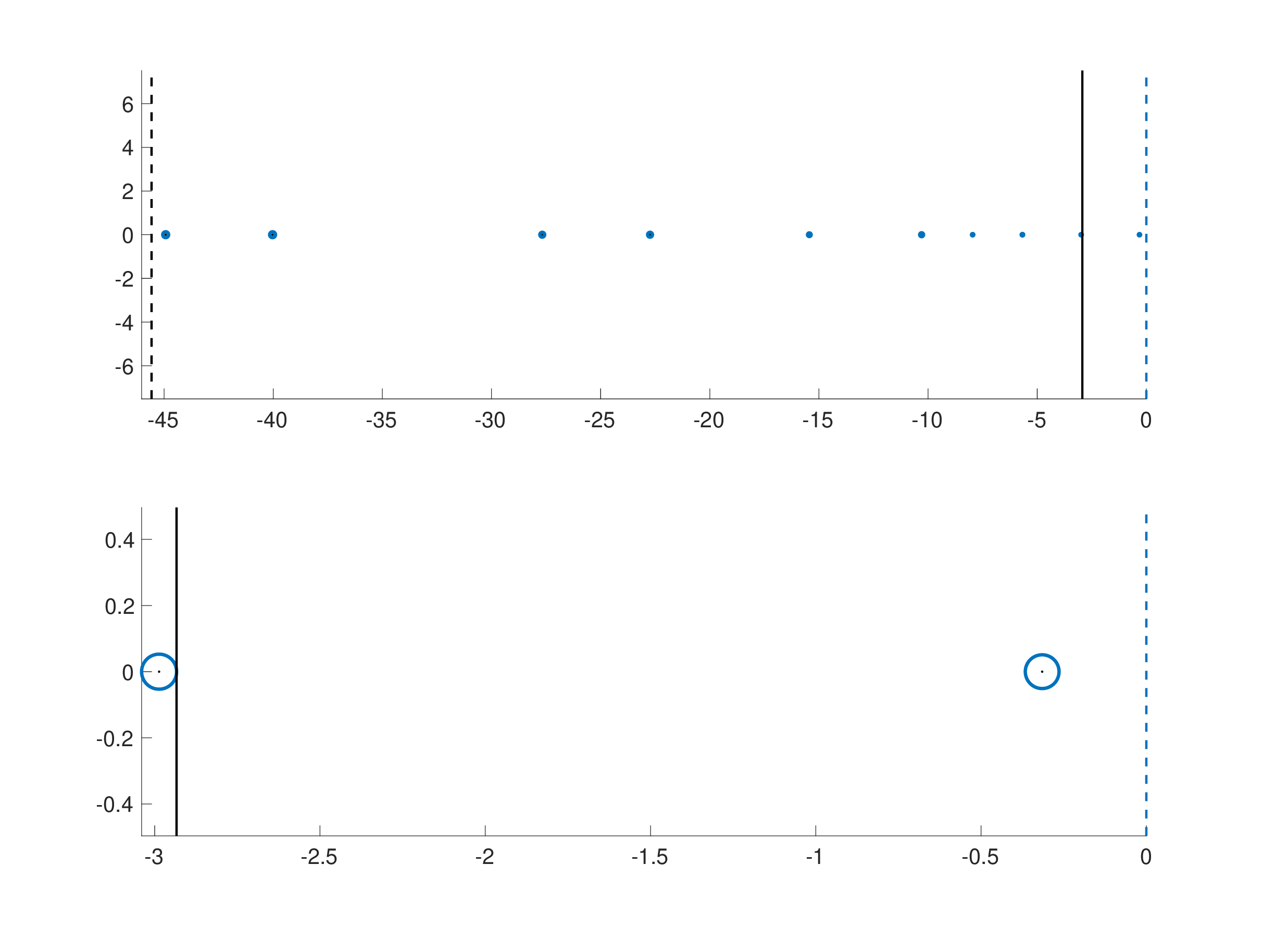}
    \put(130,600){\footnotesize $x = \overline{m}_{N,p}$}
    \put(770,600){\footnotesize $x = \mu$}
    \put(900,600){\color{myblue}\footnotesize $x = 0$}
    \put(900,270){\color{myblue}\footnotesize $x = 0$}
    \put(150,270){\footnotesize $x = \mu$}
    \put(750,195){\color{myblue}\footnotesize $\mathcal{D}(c_0,R_0)$}
    \end{overpic}
    \caption{Example of Gershgorin disks and bounds for $N=5$, $p=1.7$, $\delta = 1$, implying that all eigenvalues are stable. \\}
    \label{fig:circ1}
    \end{subfigure}
    \begin{subfigure}{\linewidth}
    \centering
    \begin{overpic}[width=0.75\linewidth, trim = 0 50 0 30, grid = false]{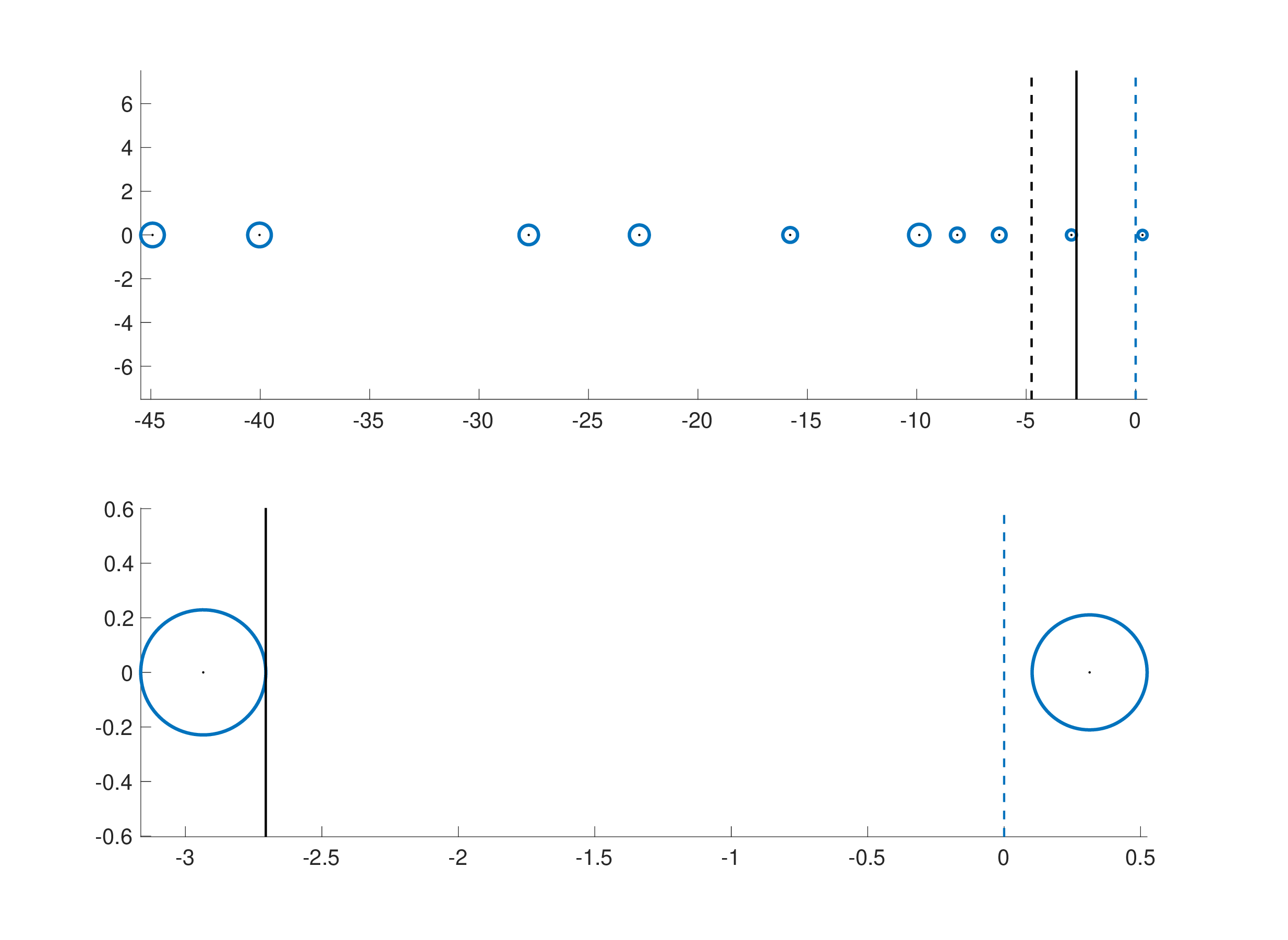}
    \put(680,600){\footnotesize $x = \overline{m}_{N,p}$}
    \put(850,570){\footnotesize $x = \mu$}
    \put(900,600){\color{myblue}\footnotesize $x = 0$}
    \put(710,260){\color{myblue}\footnotesize $x = 0$}
    \put(220,260){\footnotesize $x = \mu$}
    \put(790,225){\color{myblue}\footnotesize $\mathcal{D}(c_0,R_0)$}
    \end{overpic}
    \caption{Example of Gershgorin disks and bounds for $N=5$, $p=1.7$, $\delta = 4$, implying the existence of a single unstable eigenvalue.}
    \label{fig:circ2}
    \end{subfigure}
    \caption{Two examples showing the first $2N$ Gershgorin disks and the bound $\overline{m}_{N,p}$, illustrating the proof of Theorem~\ref{thm:spectrum}. 
    In each case, the second picture is a zoom in close to the origin. 
    Note that we intentionally took $N$ ten times smaller than in the proof of Theorem~\ref{thm:spectrum}, and $p = 1.7$, in order to get disks that are not too small and therefore easier to visualize. 
    The estimates obtained in the proof of Theorem~\ref{thm:spectrum} are actually much sharper.}
    \label{fig:circ}
\end{figure}

\begin{proof}[Proof of Theorem~\ref{thm:spectrum}] \label{proof:spectrum}
    For each $k\in\{0,\dots,999\}$, we successively consider $\delta=\left[\dfrac{4k}{1000}, \dfrac{4(k+1)}{1000}\right]$. 
    We have $q_1=q_2=1$ and take $p=2$ and $N=50$. Using the estimates from Proposition~\ref{prop:radiibounds} and Lemma~\ref{lem:last_disks}, we compute 
    $$ \mu := \max(\overline{m}_{N,p}, \{\Re(\M_{i,i})+R_i\}_{1\leq i\leq 2N-1}).$$
    We then check that $\mu < 0$ and $\Re(\M_{0,0})-R_0 > \mu$.
    This shows that $\mathcal{D}_0(\M)$ lies to the right of $\mu$, whereas all the other Gershgorin disks lie to the left of $\mu$, and by Theorem~\ref{thm:Gersh}, point $1.$ of Theorem~\ref{thm:spectrum} is proven.
    
    We now focus exclusively on $d_0(\delta)$. 
    For all $k\in\{0,\ldots,606\}$ (corresponding to all $\delta\in[0,\delta_0]$), we check that 
    $\Re(\M_{0,0})+R_0 < 0$, which proves point $2.$ of Theorem~\ref{thm:spectrum}. 
    Similarly, for all $k\in\{615,\ldots,999\}$ (corresponding to all $\delta\in[\delta_1,4]$), we check that 
    $\Re(\M_{0,0})-R_0 > 0$, which proves point $3.$ of Theorem~\ref{thm:spectrum}.
\end{proof}
The values obtained for $\mu$, and the enclosures for $d_0$ given by $d_0^-:= \Re(\M_{0,0})-R_0$ and $d_0^+:= \Re(\M_{0,0})+R_0$, are illustrated on Figure~\ref{fig:d0}. 
Examples of the $2N$ first Gershgorin disks together with the value of $\overline{m}_{N,p}$ are shown on Figure~\ref{fig:circ}.  

\section{Existence of the threshold}
\label{sec:threshold}

In order to obtain the final result of Theorem~\ref{thm:main}, we have to prove that $\delta \mapsto d_0(\delta)$ crosses $0$ only once in $(\delta_0,\delta_1)$, which corresponds to the \textbf{black} region in Figure~\ref{fig:d0} for which Theorem~\ref{thm:spectrum} does not yield precise enough information. 
To that end, we first apply a method to obtain better bounds on the first eigenvalue, based on the Newton Kantorovich Theorem \cite{Ort68}. 
Secondly, we use the same material to go further in the analysis with the implicit function theorem to finally conclude on the existence and uniqueness of the transition value $\delta^*$.

Recall that, from Corollary~\ref{cor:d0cont}, $d_0$ is a real-valued continuous function. We affirm the existence of a threshold $\delta^*$.

\begin{thm} \label{thm:threshold}
    Repeat the assumptions of Theorem~\ref{thm:spectrum}. There exits a unique $\delta^* \in (\delta_0,\delta_1)$ such that $d_0(\delta^*) = 0$. 
    Furthermore, $d_0(\delta^*) < 0$ for all  $\delta \in [\delta_0,\delta^*)$ and $d_0(\delta^*) > 0$ for all $\delta \in (\delta^*,\delta_1]$.
\end{thm}

\begin{remark}
     In fact, $\delta^* \approx 2.44456$, but we made no effort to get a tighter rigorous enclosure than $\delta^*\in(\delta_0,\delta_1)$, as this was already sufficient to prove the uniqueness of the transition.
\end{remark}

To prove Theorem~\ref{thm:threshold}, we obtain $C^1$ enclosures on the map $\delta \mapsto d_0(\delta)$ which we summarize in the following Proposition~\ref{prop:lambda'}, which then directly implies Theorem~\ref{thm:threshold}.

\begin{proposition} \label{prop:lambda'}
The function $\delta \mapsto d_0(\delta)$ is continuous on $[\delta_0,\delta_1]$ and piece-wise differentiable on $[\delta_0,\delta_1] = \bigcup_{k=607}^{614} \left[ \frac{4k}{1000}, \frac{4(k+1)}{1000} \right]$ with 
    \begin{align}
        &d_0(\delta_0) \leq-1.5\times 10^{-3} \quad \text{ and } \quad d_0(\delta_1) \geq 1.3\times 10^{-3}, \label{eq:lambda} \\
        &d_0' \left( \left[\delta_0,\delta_1 \right]\right)\subset  [0.16, 0.26]. \label{eq:lambda'}
    \end{align}
    In particular, $d_0$ is increasing on $[\delta_0,\delta_1]$.
\end{proposition}
The remainder of this section is devoted to the proof of Proposition~\ref{prop:lambda'}.


\subsection{The map \texorpdfstring{$\delta \mapsto d_0(\delta)$}{}. First properties} \label{sec:NK}
In this section, we propose a way to enclose the map $\delta \mapsto d_0(\delta)$. 
Before that, we recall that $\delta \mapsto d_0(\delta)$ is continuous, Corollary~\ref{cor:d0cont}.

Proposition~\ref{prop:lambda1} shows more than we need in Proposition~\ref{prop:lambda'}, estimates~\eqref{eq:lambda}, but it is still of interest to us to get a sharper bound on the threshold. 
In addition, the section describes a methodology for obtaining bounds on an entire curve. 
Our goal here is to remain as elementary as possible, therefore we simply use a piece-wise constant enclosure of $d_0$, combined with the implicit function theorem to get a piece-wise constant enclosure of $d_0'$. 
We point out that there exist more sophisticated and accurate methodologies to rigorously enclose curves, based on high order Chebyshev approximations but requiring a slightly more abstract framework~\cite{Bre23}, which is not necessary here.

\begin{figure}[!htbp]
    \centering
    \begin{overpic}[width=0.8\textwidth, trim = 0 100 0 0, grid = false]{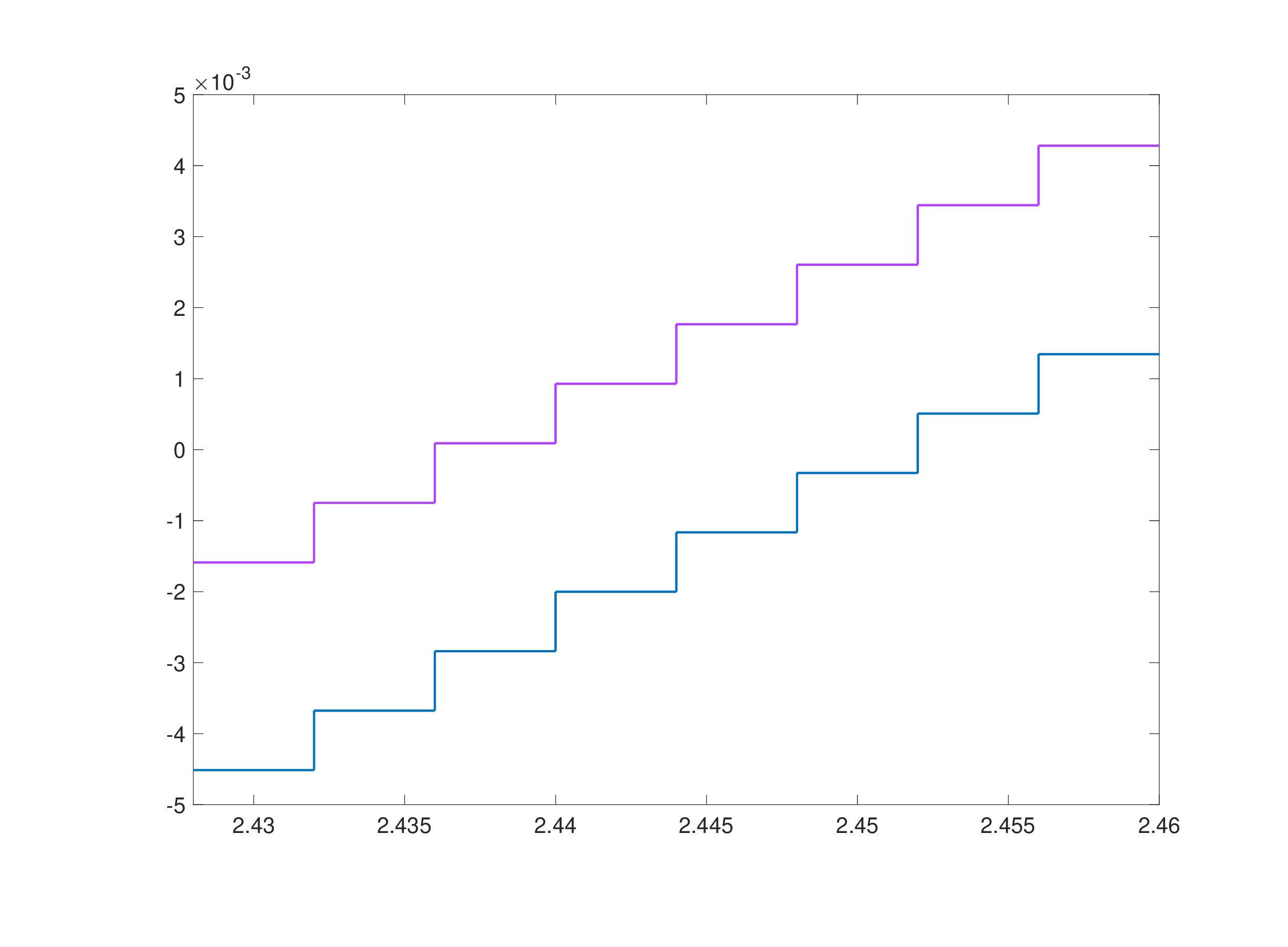}
    \put(600,0){$\delta$}
    \put(120,320){$y$}
    \put(570,200){\color{myblue}$y=\underline{d_0}(\delta)$}
    \put(470,420){\color{myviolet}$y=\overline{d_0}(\delta)$}
    \end{overpic}
    \caption{Upper- and lower-bounds on the map $\delta \mapsto d_0(\delta)$ for $\delta \in\left[\delta_0,\delta_1\right]$, with $k=607,\dots,614$.  
    These are the functions \textcolor{myblue}{$\underline{d_0}$} and \textcolor{myviolet}{$\overline{d_0}$} from Proposition~\ref{prop:lambda1}, which give sharper enclosures on $d_0$ than the ones obtained in Section~\ref{sec:Gersh} using Gershgorin disks, and shown on Figure~\ref{fig:d0}.}
    \label{fig:d0_focus}
\end{figure}
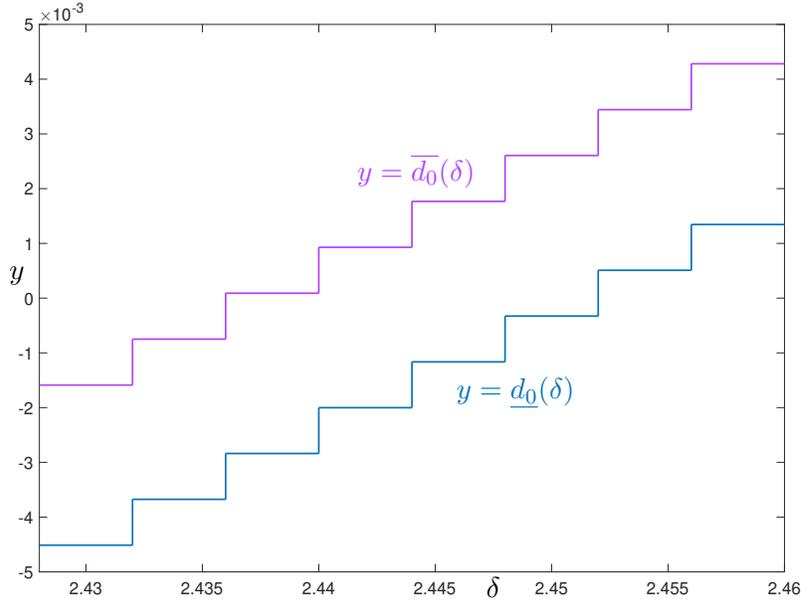

\begin{proposition} \label{prop:lambda1}
    There exist two piece-wise constant functions, called $\underline{d_0}$ and $\overline{d_0}$, defined by 
    \begin{align*}
        \underline{d_0} : [\delta_0, \delta_1] &\longrightarrow \R\\
                          \delta  &\longmapsto \underline{d_0^k}, \text{ if } \delta\in\left[\frac{4k}{1000},\frac{4(k+1)}{1000}\right),\, k = 607,\dots 614,\\
                                  &\hphantom{\longmapsto} \underline{d_0^{614}}, \text{ if } \delta=\delta_1, \\      
        \overline{d_0} : [\delta_0, \delta_1] &\longrightarrow \R\\
                          \delta  &\longmapsto \overline{d_0^k}, \text{ if } \delta\in\left[\frac{4k}{1000},\frac{4(k+1)}{1000}\right),\, k = 607,\dots 614,\\
                          &\hphantom{\longmapsto} \overline{d_0^{614}}, \text{ if } \delta=\delta_1,
    \end{align*}
    
    \begin{table}[!htbp]
        \centering
        \begin{equation*}
        \begin{array}{c|c|c|c|c|c|c|c|c}
            k & 607 & 608 & 609 & 610 & 611 & 612 & 613 & 614 \\ \hline
            \underline{d_0^k} \times 10^3 & -4.6 & -3.7 & -2.9 & -2.1 & -1.2 & -0.4 & 0.5 & 1.3 \\[0.5em] 
            \overline{d_0^k} \times 10^3 & -1.5 & -0.7 & 0.1 & 1.0 & 1.8 & 2.7 & 3.5 & 4.3 \\ 
        \end{array}
        \end{equation*}
        \caption{Description of the two piece-wise constant functions $\underline{d_0}$ and $\overline{d_0}$ of Proposition~\ref{prop:lambda1}.}
        \label{tab:enclose}
    \end{table}
    such that
    \begin{equation}
        \forall\, \delta \in [\delta_0,\delta_1], \ \underline{d_0}(\delta) \leq d_0(\delta) \leq \overline{d_0}(\delta).
        \label{eq:enclose}
    \end{equation}
    The functions $\underline{d_0}$ and $\overline{d_0}$ are depicted in Figure~\ref{fig:d0_focus} and  defined in Table~\ref{tab:enclose}. 
\end{proposition}

\begin{remark}
    From Proposition~\ref{prop:lambda1}, we have the inequalities \eqref{eq:lambda} of Proposition~\ref{prop:lambda'} satisfied.
    Indeed, $d_0(\delta_0)\leq \overline{d_0}(\delta_0) = \overline{d_0^{607}}=-1.5\times10^{-3}$ and  $d_0(\delta_1) \geq \underline{d_0}(\delta_1) = \underline{d_0^{614}}= 1.3\times10^{-3}$.
\end{remark}

The remainder of Section~\ref{sec:NK} is devoted to the proof of Proposition~\ref{prop:lambda1}. 
We first fix some notation and introduce in Section~\ref{subsec:new_pov} a zero-finding problem whose zeros are eigenpairs of $M$. 
In order to study zeros of the problem, we use a Newton-like fixed point operator in Section~\ref{sec:thmNK}, with a suitable approximate inverse built in Section~\ref{sec:A}. 
We then derive quantitative bounds to study this fixed point operator in Section~\ref{sec:boundsYZ1Z2}, and use them in Section~\ref{sec:prooflambda1} to prove Proposition~\ref{prop:lambda1}.

\subsubsection{Change of basis and definitions: a new point of view} \label{subsec:new_pov}

We get back to system~\eqref{eq:system_infinite}. 
We denote $u(t) = (u_k(t))_{k\in\N}, \ v(t) = (v_k(t))_{k\in\N}$ the sequences of Fourier modes of the functions $u(\cdot,t),\ v(\cdot,t)$, for any $t\in\R_+$. 
Thus, the dynamics in \eqref{eq:system_infinite} can be written as follow
\begin{align}
    \begin{cases}
        u' = \vartheta\Delta u + a u + b v, \\
        v' = \Delta v + cu + dv + \delta B u, \\
    \end{cases}
\end{align}
where $\Delta : (u_k) \mapsto (-(\frac{k\pi}{l})^2 u_k)$ and $B = (B_{i,j})_{(i,j)\in\N^2}$ as in \eqref{eq:B1}.

Then we denote $M$ the linear operation on $\begin{pmatrix} u\\v \end{pmatrix}$, $M = \begin{pmatrix}
    \vartheta\Delta + aI & bI \\
    cI + \delta B & \Delta + dI 
\end{pmatrix}$.
Note that $M$ depends on $\delta$ but so as to lighten the notations we do not include it.

\begin{remark}
We note that, although we kept the same letter $M$ for simplicity, this infinite matrix is not exactly the same $M$ as the one in Section~\ref{sec:Gersh}. 
The only difference is a change of basis given by a permutation: in Section~\ref{sec:Gersh} we wrote $M$ assuming the Fourier coefficients of $u$ and $v$ were ordered $(u_0, v_0, u_1, v_1, \dots )$, but in this section it will be more convenient to split the two components and have $(u_0, u_1, \dots ; v_0, v_1, \dots) = (u,v)$. 
Of course one could in principle stick with the same convention for both sections, but we felt that each of them made their respective section easier to follow.
\end{remark}

We are looking for $d_0$, but we first present a technique that allows to very precisely enclose any single eigenvalue $\lambda$ of $M$, but without telling us \emph{which} eigenvalue we are enclosing. 
However, note that Theorem~\ref{thm:spectrum} already provides us with some control on $d_0$. In particular, we have a threshold $\mu = \mu(\delta)$ (obtained explicitly in the proof of Theorem~\ref{thm:spectrum}) such that, for all $\delta\in[0,4]$, $d_0$ is the only eigenvalue of $M$ whose real part is larger than $\mu$.  
Therefore, once a very accurate enclosure of an eigenvalue $\lambda$ of $M$ is rigorously obtained, we can a posteriori prove that $\lambda$ indeed corresponds to $d_0$ by checking that $\Re(\lambda)>\mu$.

Let us now consider the eigenvalue-eigenvector problem
\begin{equation} \label{eq:eig_prob}
    M\begin{pmatrix} u\\ v\end{pmatrix} = \lambda \begin{pmatrix} u\\ v\end{pmatrix},
\end{equation}
where $(\lambda,u,v)\in\C\times\C^\N\times\C^\N$ are unknowns. 
In order to get an isolated solution, we add a normalization condition to that system. 
Therefore, we first compute numerically a unitary finite approximate eigenvector $(\tilde{u},\tilde{v})$ related to the eigenvalue $d_0$,
and we then search for zeros of the following functional:
\begin{align}
    F(\lambda,u,v) = \left( \begin{array}{c}
    \begin{pmatrix} \tilde{u} \\ \tilde{v} \end{pmatrix} \!\cdot\! \begin{pmatrix} u \\ v \end{pmatrix} - 1 \\
    M \begin{pmatrix} u \\ v \end{pmatrix} - \lambda \begin{pmatrix} u \\ v \end{pmatrix}
\end{array}\right).
\end{align}

If $(\lambda,u,v)$ is a zero of $F$, then $\lambda$ is an eigenvalue of $M$.
By rewriting the search for $\lambda$ as a zero-finding problem, we will be able to leverage by now standard and powerful computer-assisted techniques~\cite{BerLes15} to get a really tight and rigorous enclosure of $\lambda$. 
This zero-finding setting also allows to study parameter dependency, therefore we will naturally be able to use the implicit function theorem to then study the derivative of $\lambda$ with respect to $\delta$. 
Before presenting in more details these computer-assisted techniques, we introduce some notations.

\begin{definition}
   Let $\alpha \in\R$, let $\ell^1_\alpha$ the set of sequences $u = (u_k)_{k\in\N}$ such that:
$$\| u\|_\alpha := |u_0| + 2\sum_{k=1}^{+\infty} |u_k| k^\alpha < +\infty.$$
Then $\ell^1_\alpha$ is a Banach space. 
Furthermore, if $\alpha \geq 0$, the weights are sub-multiplicative and it brings to $\ell^1_\alpha$ a structure of a Banach algebra with the discrete convolution, see \cite{Les18}.
\end{definition}

\begin{definition}
    We denote $\mathcal{X}_\alpha = \C\times \ell^1_\alpha\times\ell^1_\alpha$. We define the following norm: $$X = (\lambda,u,v) \in \mathcal{X}_\alpha,\ \|X\|_{\mathcal{X}_\alpha} = |\lambda| + \|u\|_\alpha + \|v\|_\alpha.$$
\end{definition}

\begin{remark}
    $\mathcal{X}_\alpha$ is isomorphic to $\ell^1$, the $\mathcal{X}_\alpha$-norm is a weighted $\ell^1$-norm. $\mathcal{X}_\alpha$ is a Banach space.
\end{remark}

\begin{definition}
\label{def:blockop}
    The operator norm induced on $\mathcal{X}_\alpha$ can be expressed as follows.
    Let  $l\in\C, \ l_1,\tilde{l}_1, l_2, \tilde{l}_2 \in \ell^1_\alpha$ and $L^{11},L^{12},L^{21}, L^{22} \in \mathcal{L}(\ell^1_\alpha)$. 
    Write
    \begin{align*}
    L =  
    \left(\begin{array}{c|c|c}
    l & \tilde{l}_1^* & \tilde{l}_2^* \\ \hline
    l_1 & L^{11} &  L^{12} \\ \hline
    l_2 &  L^{21} &  L^{22}
    \end{array}\right)
    \in\mathcal{L}(\mathcal{X}_\alpha),
    \end{align*}
    where $\tilde{l}_1^*$ and $\tilde{l}_2^*$ are the adjoints of $\tilde{l}_1$ and $\tilde{l}_2$, which are linear forms on $\ell^1_\alpha$.
    
    The operator norm of $L$ induced by $\|\cdot\|_{\mathcal{X}_\alpha}$ is 
    \begin{align}
        \| L \|_{\mathcal{L}(\mathcal{X}_\alpha)} &= \sup_{X\in\mathcal{X}_\alpha, \|X\|_{\mathcal{X}_\alpha} > 0} \dfrac{\|LX\|_{\mathcal{X}_\alpha}}{\|X\|_{\mathcal{X}_\alpha}}, \notag \\
        &= \max \left\{\left\|\begin{array}{c} l \\ \hline l_1 \\ \hline l_2 \end{array} \right\|_{\mathcal{X}_\alpha}; \sup_{j\geq 0} \frac{1}{\max(1,2j^\alpha)} \left\|\begin{array}{c} (\tilde{l_1})_j \\ \hline (L^{11})_j \\ \hline (L^{21})_j \end{array} \right\|_{\mathcal{X}_\alpha}; \sup_{j\geq 0} \frac{1}{\max(1,2j^\alpha)} \left\|\begin{array}{c} (\tilde{l_2})_j \\ \hline (L^{12})_j \\ \hline (L^{22})_j \end{array} \right\|_{\mathcal{X}_\alpha} \right\}, \label{def:op_norm}
    \end{align}
    where $(\, \cdot\, )_j$ denotes the $j^{th}$ component of an element of $\ell^1_\alpha$ (it is then an element of $\C$) or of $\mathcal{L}(\ell^1_\alpha)$ (it is then an element of $\ell^1_{\alpha}$).
\end{definition}

\begin{remark}
To simplify the reading of this formula we introduce the norm of each column block,
\begin{align}
    C^{\eta}_\alpha(L) := 
        \left\{\begin{array}{ll}
        \left\|\begin{array}{c} l \\ \hline l_1 \\ \hline l_2 \end{array} \right\|_{\mathcal{X}_\alpha}, & \text{ if } \eta = 0 \\
        \sup_{j\geq 0} \frac{1}{\max(1,2j^\alpha)} \left\|\begin{array}{c} (\tilde{l_1})_j \\ \hline (L^{11})_j \\ \hline (L^{21})_j \end{array} \right\|_{\mathcal{X}_\alpha}, & \text{ if } \eta = 1 \\
        \sup_{j\geq 0} \frac{1}{\max(1,2j^\alpha)} \left\|\begin{array}{c} (\tilde{l_2})_j \\ \hline (L^{12})_j \\ \hline (L^{22})_j \end{array} \right\|_{\mathcal{X}_\alpha}, & \text{ if } \eta = 2 \end{array} \right.
\end{align}
Hence, we have $\| L \|_{\mathcal{L}(\mathcal{X}_\alpha)} = \max\limits_{\eta = 0,1,2} C^{\eta}_\alpha(L)$.
\end{remark}

\begin{remark}
    Based on~\eqref{def:op_norm}, as soon as we can either compute or estimate the norms of all the columns of a linear operator $L$, then we can compute or at least explicitly estimate its operator norm. 
    This will be important for the main results of Section~\ref{sec:threshold}.
\end{remark}

\subsubsection{Zero of \texorpdfstring{$F$}{}}
\label{sec:thmNK}

Our goal is now to get a precise and rigorous description of a zero of the map $F$ defined from $\mathcal{X}_\alpha$ into $\mathcal{X}_{\alpha-2}$.

Starting from an approximate zero $\bar{X}=(\bar{\lambda},\bar{u},\bar{v})$ of $F$ computed numerically, our goal will be to establish the existence of a nearby exact zero of $F$, and to give an explicit error bound. 
We will use a version of the Newton-Kantorovich theorem, which requires first looking at the Fréchet derivative of $F$.
Following the notations of Definition~\ref{def:blockop}, the Fréchet derivative is
\begin{align}
DF(X) = \left(\begin{array}{c|c|c}
0 & \tilde{u}^* & \tilde{v}^* \\ \hline
-u & \vartheta \Delta + (a-\lambda)I_{\ell^1} & bI_{\ell^1} \\ \hline
-v & cI_{\ell^1} +\delta B & \Delta + (d-\lambda)I_{\ell^1} 
\end{array}\right).
\end{align}
$DF(X)$ is a linear operator from $\mathcal{X}_\alpha$ to  $\mathcal{X}_{\alpha-2}$. 
We treat each block as an operator: $0$ as the multiplication by $0$ from $\C$ into $\C$; $\tilde{u}^*$ (resp. $\tilde{v}^*$) as the linear form corresponding to the dual element of $\tilde{u}$ (resp. $\tilde{v}$) from $\ell^1_\alpha$ into $\C$; $-u$ (resp. $-v$) as the vector multiplication by $-u$ (resp. $-v$) from $\C$ into $\ell^1_\alpha$; and the last blocks are operators from $\ell^1_\alpha$ into $\ell^1_{\alpha-2}$.

In order to rigorously establish the existence of a zero of $F$ near $\bar{X}$, we consider an injective approximate inverse $A \in \mathcal{L}(\mathcal{X}_{\alpha-2},\mathcal{X}_\alpha)$ of $DF(X)$ (see Definition~\ref{def:A} for a more precise description of $A$).
The following theorem, which is common in computer-assisted proofs in nonlinear analysis~\cite{Plu92,Yam98,DayLesMis07,AriKoc10}, provides sufficient conditions to prove that the operator
\begin{align*}
    \begin{cases}
    \mathcal{X}_\alpha \to \mathcal{X}_\alpha, \\
    X \mapsto X-AF(X),
\end{cases}
\end{align*}
is a contraction on a small and explicit neighborhood of $\bar{X}$, which proves the existence of a zero of $F$ near $\bar{X}$, and gives explicit error bounds between that zero and $\bar{X}$.

\begin{thm} \label{thm:NK}
Let $\alpha \geq 0$, let $\bar{X} \in \mathcal{X}_\alpha$, and let $A$ be as in~\eqref{def:A}. Let $Y, Z_1, Z_2$ positive such that:
\begin{align}
    \|AF(\bar{X})\|_{\mathcal{X}_\alpha} &\leq Y, \label{hyp:NK-Y}\\
    \|I_{\mathcal{X}} - ADF(\bar{X})\|_{\mathcal{L}(\mathcal{X}_\alpha)} &\leq Z_1, \label{hyp:NK-Z1}\\ 
    \| A(DF(\bar{X}) - DF(X)) \|_{\mathcal{L}(\mathcal{X}_\alpha)} &\leq Z_2 \|\bar{X}-X\|_{\mathcal{X}_\alpha}, \ \forall X \in \mathcal{X}_\alpha .\label{hyp:NK-Z2}
\end{align}
If we have $(1-Z_1)^2 - 2Z_2Y > 0$ and $Z_1<1$,
we denote $r_{\min} = \frac{(1-Z_1) - \sqrt{(1-Z_1)^2 - 2Z_2Y}}{Z_2}$ and $r_{\max} = \frac{1-Z_1}{Z_2}$.
Then, for any $r \in [r_{\min},r_{\max})$, there exists a unique $\widehat{X} \in \mathcal{B}_{\mathcal{X}_\alpha}(\bar{X},r)$ such that $F(\widehat{X}) = 0$.
\end{thm}

\begin{remark}
    $F$ is a quadratic functional and its second derivative norm is equal to $1$. 
    The hypotheses~\eqref{hyp:NK-Z2} is therefore satisfied as soon as $\|A\|_{\mathcal{L}(\mathcal{X}_\alpha)} \leq Z_2$. 
    This theorem can easily be generalized to situation where $F$ has higher order nonlinearities, but this is not needed here.
\end{remark}

\begin{remark}
    Assume there exists $Y,Z_1,Z_2$ positive and $\bar{X}=(\bar{\lambda},\bar{u},\bar{v})$ an approximate zero of $F$, such that all hypotheses of Theorem~\ref{thm:NK} are satisfied.
    Let $\widehat{X}=(\widehat{\lambda},\widehat{u},\widehat{v})$ be the theoretical zero given by the conclusion of Theorem~\ref{thm:NK}.
    We have $||\bar{X}-\widehat{X}||_{\mathcal{X}_\alpha} \leq r_{\min}$, thus $|\bar{\lambda} - \widehat{\lambda}| \leq r_{\min}$. In other words, we get an \emph{a posteriori} error bound on the approximate eigenvalue $\bar{\lambda}$, which provides us with an explicit and guaranteed enclosure $[\bar{\lambda}-r_{\min},\bar{\lambda}+r_{\min}]$ of the exact eigenvalue $\hat\lambda$.
\end{remark}

\subsubsection{Formalism and construction of \texorpdfstring{$A$}{}}
\label{sec:A}

To complete our demonstration, we need to show how we get the hypotheses of Theorem~\ref{thm:NK} for any $\delta$. 
The statement of the hypotheses of Theorem~\ref{thm:NK} is based on the operator $A$. 
To build it, we use again the truncation operator. 
Firstly, we have to explain how it interacts with the new point of view, introduced in Section~\ref{subsec:new_pov}.

\begin{definition}
    Let $z\in\ell^1_\alpha$, we denote $\Pi^N_{\ell^1} z$ the element of $\ell^1_\alpha$ such that $$(\Pi^N_{\ell^1}z)_k = \begin{cases} z_k, \ & k \ < N, \\ 0, \ &k \geq N. \end{cases}$$
    We identify $\Pi^N_{\ell^1}z$ with an element of $\C^N$.
    $\Pi^N_{\ell^1}$ is an operator on $\ell^1_\alpha$ (for any $\alpha$). The operator $\Pi^N_{\ell^1}$ can be interpreted as an infinite matrix,
    $(\Pi^N_{\ell^1})_{i,j} = \begin{cases}
        1, &\text{ if } i=j \text{ and } i < N, \\
        0, &\text{ otherwise.}
    \end{cases}$
\end{definition}
We extend this notation to any element of $\mathcal{X}_\alpha$.
\begin{definition}\label{def:pi_X}
    For any $Z = (\mu, y, z) \in \mathcal{X}_\alpha$, 
    \begin{align*}
        \Pi^N_{\mathcal{X}}Z &= (\mu, \Pi^N_{\ell^1}y, \Pi^N_{\ell^1}z) \in \C\times\C^N\times\C^N.
    \end{align*}
$\Pi^N_{\mathcal{X}}$ is an operator on $\mathcal{X}_\alpha$ (for any $\alpha$), being interpreted as 
$3\times 3$ block matrix:
$\Pi^N_{\mathcal{X}} = \left( \begin{array}{c|c|c}
    1 & 0 & 0 \\ \hline
    0 & \Pi^N_{\ell^1} & 0 \\ \hline
    0 & 0 & \Pi^N_{\ell^1}
\end{array}\right)$.
\end{definition}

\begin{remark}
    After the definitions of the projectors, we denote for any $u\in\ell^1_\alpha$ and $X\in\mathcal{X}_\alpha$ finite vectors by writing $u\in\Pi^N_{\ell^1} \ell^1_\alpha$ and $X\in\Pi^N_\mathcal{X} \mathcal{X}_\alpha$.
\end{remark}

\begin{remark}
  In the $3\times3$ block matrix of Definition~\ref{def:pi_X}, the zero operators have different domains and range, but to simplify notations we make no distinctions between them.
\end{remark}
We give below the examples of truncation for $B$ and $\Delta$.
\begin{example}
Let $\Delta : \begin{cases}
     \ell^1_\alpha \to \ell^1_{\alpha-2}  \\
     (u_k) \mapsto (-(\frac{k\pi}{l})^2 u_k)
\end{cases}$, which is a diagonal operator. We have 
\begin{align*} \Pi^N_{\ell^1}\Delta &: \begin{cases}
     \ell^1_\alpha \to \ell^1_{\alpha-2}  \\
     (u_k) \mapsto \begin{cases}
         -(\frac{k\pi}{l})^2 u_k, \ k<N \\
         0,\ k\geq N.
     \end{cases}
\end{cases} \\
(I_{\ell^1}-\Pi^N_{\ell^1})\Delta &: \begin{cases}
     \ell^1_\alpha \to \ell^1_{\alpha-2}  \\
     (u_k) \mapsto \begin{cases}
         0, \ k<N \\
         -(\frac{k\pi}{l})^2 u_k,\ k\geq N.
     \end{cases}
\end{cases}
\end{align*}
\end{example}
We see here that $\Pi^N_{\ell^1}$ and $\Delta$ commute, since $\Delta$ is diagonal.

\begin{example}
Let $B : \begin{cases}
     \ell^1_\alpha \to \ell^1_\beta  \\
     (u_k) \mapsto (\sum_{j=0}^{+\infty}B_{k,j}u_j)_{k\in\N}
\end{cases}$, where $B$ satisfies \eqref{hyp:B} with $\alpha \geq -q_2, \ \beta < q_1-1$. We have 
\begin{align*}
    \Pi^N_{\ell^1} B &: \begin{cases}
     \ell^1_\alpha \to \ell^1_\beta  \\
     (u_k) \mapsto \begin{cases}
         \sum_{j=0}^{+\infty}B_{k,j}u_j,\ k < N, \\
         0,\ k \geq N.
     \end{cases} 
\end{cases} \\
(I_{\ell^1}-\Pi^N_{\ell^1}) B &: \begin{cases}
     \ell^1_\alpha \to \ell^1_\beta  \\
     (u_k) \mapsto \begin{cases}
        0,\ k < N, \\
        \sum_{j=0}^{+\infty}B_{k,j}u_j,\ k \geq N.
     \end{cases} 
\end{cases}
\end{align*}
\end{example}
Be careful, here $B$ and $\Pi^N_{\ell^1}$ do not commute.

The property “to commute with $\Pi^N_{\ell^1}$" means that the operator do not mix the $N$ first Fourier modes and the tail of any sequence. 
This is a significant property to keep in mind in the further calculation. With a matrix point of view, an operator $ L$ on $\ell^1_\alpha$ commutes with $\Pi^N_{\ell^1}$ if and only if $L$ is a $2\times 2$ block diagonal matrix:
\begin{equation*}
     L \ \text{commutes with} \ \Pi^N_{\ell^1} \iff  L = \left(\begin{array}{c|c}
     L_N & \\ \hline
     &  L_\infty
     \end{array}\right)
     \begin{array}{l}
          \leftarrow \text{operates on the first $N$ modes} \\
          \leftarrow \text{operates on the tail of sequences} 
     \end{array}
\end{equation*}
The following lemma gives the algebraic translation.

\begin{lem} \label{lem:commute} Let $\alpha,\beta \in \R$.
Let $ L\in \mathcal{L}(\ell^1_\alpha,\ell^1_\beta)$ and $u,v\in \ell^1_\alpha$, we have
\begin{align*} 
    \Pi^N_{\ell^1}L =L\Pi^N_{\ell^1} 
    &\iff \Pi^N_{\ell^1} L (I_{\ell^1}-\Pi^N_{\ell^1}) = 0 \text{ and }(I_{\ell^1}-\Pi^N_{\ell^1}) L \Pi^N_{\ell^1} = 0 \\
    &\iff (I_{\ell^1}-\Pi^N_{\ell^1})L = L (I_{\ell^1}-\Pi^N_{\ell^1}).
\end{align*}
\end{lem}

We use the formalism introduce in the beginning of section~\ref{sec:A} to show how we built $A$, that appears in Theorem~\ref{thm:NK}. 
The idea behind Definition~\ref{def:A} is to get, in finite dimension, a good approximation of $DF(X)^{-1}$ and to only  keep for the complement the main part of the inverse of $(I_{\mathcal{X}} - \Pi^N_{\mathcal{X}})DF(X)(I_{\mathcal{X}} - \Pi^N_{\mathcal{X}})$, which is “easy" to deal with (e.g. diagonal).

\begin{definition} \label{def:A}
Let $\alpha \in [-q_2, q_1+1)$, with $q_1,q_2$ the parameters from \eqref{hyp:B}.
Let $\bar{X} = (\bar{\lambda},\bar{u},\bar{v}) \in \Pi^N_\mathcal{X}\mathcal{X}_\alpha$, and $\tilde{u},\tilde{v}  \in \Pi^N_{\ell^1}\ell^1_\alpha$.

We choose $A\in\mathcal{L}(\mathcal{X}_{\alpha-2},\mathcal{X}_\alpha)$, in the following way:
\begin{align}
    &\Pi^N_{\mathcal{X}} A \Pi^N_{\mathcal{X}}= A_N \\
    \intertext{where}
    &A_N = \left(\begin{array}{c|c|c}
        a_0 & \tilde{a}_1^* & \tilde{a}_2^* \\ \hline
        a_1 & A^{11} & A^{12} \\ \hline
        a_2 & A^{21} & A^{22}
\end{array}\right) \approx \left(\Pi^N_{\mathcal{X}} DF(\bar{X}) \Pi^N_{\mathcal{X}}\right)^{-1},  \label{eq:A} \\
&A - A_N = \left(\begin{array}{c|c|c}
        0 & 0 & 0 \\ \hline
        0 &(I_{\ell^1}-\Pi^N_{\ell^1}) (\vartheta\Delta + (a-\bar{\lambda})I_{\ell^1})^{-1} & 0 \\ \hline
        0 & 0 & (I_{\ell^1}-\Pi^N_{\ell^1})(\Delta + (d-\bar{\lambda})I_{\ell^1})^{-1}\notag
\end{array}\right).
\end{align}
Here $A^{ij}$ are finite matrices of size $N\times N$, hence $A_N$ is of size $(2N+1)\times (2N+1)$. Similarly, $\Pi^N_{\mathcal{X}} DF(\bar{X}) \Pi^N_{\mathcal{X}}$ can be interpreted as a $(2N+1)\times (2N+1)$ matrix, hence $A_N$ can simply be obtained numerically by computing an approximate inverse.
\end{definition}

\begin{remark}
    The reason why $\alpha$ must belong to $[-q_2, q_1+1)$ is to correctly define $DF(X)$ from $\ell^1_\alpha$ to $\ell^1_{\alpha-2}$.
\end{remark}

\subsubsection{Derivation of the bounds \texorpdfstring{$Y, Z_1, Z_2$}{}}
\label{sec:boundsYZ1Z2}

Now that $A$ has been fixed, we are ready to derive computable estimates $Y$, $Z_1$ and $Z_2$ satisfying assumptions~\eqref{hyp:NK-Y}-\eqref{hyp:NK-Z2} of Theorem~\ref{thm:NK}.

\begin{proposition} \label{prop:bounds}
    Let $a, b, c, d, \vartheta$ be fixed as in Theorem~\ref{thm:main}, and $l = |\Omega|$. 
    Let $\delta\in\R$ and $B$ satisfying \eqref{hyp:B}. Let $\alpha \in [-q_2,1+q_1)$, with $\alpha \geq 0$ and $\bar{X} = (\bar{\lambda},\bar{u},\bar{v}) \in \Pi^N_\mathcal{X}\mathcal{X}_\alpha$, and $\tilde{u},\tilde{v} \in \Pi^N_{\ell^1}\ell^1_\alpha$, with $N\in\N$ such that $N > 1 + \frac{l}{\pi}\sqrt{|d-\bar{\lambda}|}.$
    
    We denote 
    \begin{equation*}
        E(\alpha,N) = \left(\dfrac{l}{\pi}\right)^2\dfrac{(N-1-\nu)^{\alpha-(1+q_1)}}{(1+q_1)-\alpha}, \text{where } \nu = \dfrac{l}{\pi}\sqrt{\max(0,d-\bar{\lambda})}.
    \end{equation*}
    Furthermore, denote 
    \begin{equation*}
        R(\theta,x,N) = \dfrac{1}{| -\theta\left(\frac{N\pi}{l}\right)^2 + (x-\bar{\lambda})|}, \text{ where } \theta \in \{ 1, \ \vartheta\} \text{ and } x \in \{a,\ d\}.
    \end{equation*} 
    Let $\chi_q = (1,1,\frac{1}{2^q},\dots,\frac{1}{(N-1)^q})^* \in \C^N$ for $q \in \{q_1,q_2\}$. 
    In the sequel, $|\cdot|$ applied to a vector or a matrix has to be understood element-wise.
The quantities
\begin{align*}
    Y &= \|A_N\Pi^N_\mathcal{X}F(\bar{X})\|_{\mathcal{X}_\alpha} + 2C|\delta|E(\alpha,N)\left(|\bar{u}|\cdot\chi_{q_2}\right), \\
    Z_1 &= \max\Big\{ \|\Pi^N_{\mathcal{X}} - \Pi^N_{\mathcal{X}}A DF(\bar{X})\Pi^N_\mathcal{X}\|_{\mathcal{L}(\mathcal{X}_\alpha)}+2C|\delta| E(\alpha,N); \\
    &\qquad\qquad    \dfrac{C|\delta|}{2N^{\alpha+q_2}} \left( |\tilde{a}_2|\cdot\chi_{q_1} + \left\| |A^{12}|\chi_{q_1} \right\|_{\ell^1_\alpha} + \left\| |A^{22}|\chi_{q_1} \right\|_{\ell^1_\alpha} \right)  + |c| R(1,d,N) + \frac{C|\delta|}{N^{\alpha+q_2}} E(\alpha,N); \\
    &\qquad\qquad |b| R(\vartheta,a,N) \Big\}, \\
    Z_2 &= \max \left\{ \|A_N\|_{\mathcal{L}(\mathcal{X}_\alpha)}; R(\vartheta,a,N);  R(1,d,N) \right\},
\end{align*}
satisfy \eqref{hyp:NK-Y}, \eqref{hyp:NK-Z1} and \eqref{hyp:NK-Z2} respectively.

\end{proposition}

\begin{proof} We get each of the three bounds separately.

\begin{proof*}{Derivation of $Y$.}
$Y$ is deduced from $\|AF(\bar{X})\|_{\mathcal{X}_\alpha}$.
Let us bound $\|AF(\bar{X})\|_{\mathcal{X}_\alpha}$.

We have
\begin{align*}
    &F(\bar{X}) = \left(\begin{array}{c}
    \bar{u}\cdot \tilde{u} + \bar{v}\cdot \tilde{v} - 1  \\ \hline
    \vartheta \Delta \bar{u} + (a-\bar{\lambda})\bar{u} + b\bar{v} \\ \hline
    c\bar{u} +\delta B\bar{u} + \Delta \bar{v} + (d-\bar{\lambda})\bar{v} 
\end{array}\right),
\end{align*}
and then
\begin{align*}
\Pi^N_\mathcal{X}F(\bar{X}) &= 
\left(\begin{array}{c}
    \Pi^N_{\ell^1}\bar{u}\cdot \Pi^N_{\ell^1}\tilde{u} + \Pi^N_{\ell^1}\bar{v}\cdot \Pi^N_{\ell^1}\tilde{v} - 1  \\ \hline
    \vartheta \Delta \Pi^N_{\ell^1}\bar{u} + (a-\bar{\lambda})\Pi^N_{\ell^1}\bar{u} + b\Pi^N_{\ell^1}\bar{v} \\ \hline
    c\Pi^N_{\ell^1}\bar{u} +\delta \Pi^N_{\ell^1}B\Pi^N_{\ell^1}\bar{u} + \Delta\Pi^N_{\ell^1} \bar{v} + (d-\bar{\lambda})\Pi^N_{\ell^1}\bar{v} 
\end{array}\right), 
\end{align*}
but $\bar{X} \in\Pi^N_\mathcal{X}\mathcal{X}_\alpha$ and $\tilde{u},\tilde{v} \in \Pi^N_{\ell^1} \ell^1_\alpha$, so
\begin{align*}
\Pi^N_\mathcal{X}F(\bar{X}) &=
\left(\begin{array}{c}
    \bar{u}\cdot\tilde{u} + \bar{v}\cdot\tilde{v} - 1  \\ \hline
    \vartheta \Delta \bar{u} + (a-\bar{\lambda})\bar{u} + b\bar{v} \\ \hline
    c\bar{u} +\delta \Pi^N_{\ell^1}B\bar{u} + \Delta \bar{v} + (d-\bar{\lambda})\bar{v}
\end{array}\right), \\
(I_\mathcal{X} - \Pi^N_\mathcal{X})F(\bar{X}) &= \left(\begin{array}{c}
    0  \\ \hline
    0 \\ \hline
    \delta (I_{\ell^1}-\Pi^N_{\ell^1})B\Pi^N_{\ell^1}\bar{u}
\end{array}\right).
\end{align*}

Since $A$ commutes with $\Pi^N_{\mathcal{X}}$, see \eqref{eq:A}, we have with Lemma~\ref{lem:commute} $$A\Pi^N_{\mathcal{X}}F(\bar{X}) = A\Pi^N_{\mathcal{X}}\Pi^N_{\mathcal{X}}F(\bar{X}) =  \Pi^N_{\mathcal{X}}A\Pi^N_{\mathcal{X}}F(\bar{X}) = A_N\Pi^N_\mathcal{X}F(\bar{X}), $$
which is a finite vector, and can therefore be computed fully with the computer.
Similarly, 
$$A(I_{\mathcal{X}}-\Pi^N_{\mathcal{X}}) F(\bar{X}) = (I_{\mathcal{X}}-\Pi^N_{\mathcal{X}}) A(I_{\mathcal{X}}-\Pi^N_{\mathcal{X}}) F(\bar{X}) = (A-A_N)(I_{\mathcal{X}}-\Pi^N_{\mathcal{X}}) F(\bar{X}).$$
Then, we have
\begin{align*}
    \|AF(\bar{X})\|_{\mathcal{X}_\alpha} &= \| A\Pi^N_{\mathcal{X}} F(\bar{X}) + A (I_{\mathcal{X}}-\Pi^N_{\mathcal{X}}) F(\bar{X})\|_{\mathcal{X}_\alpha} \\
    &=\| A_N\Pi^N_{\mathcal{X}} F(\bar{X}) + (A - A_N) (I_{\mathcal{X}}-\Pi^N_{\mathcal{X}}) F(\bar{X})\|_{\mathcal{X}_\alpha} \\
    &\leq \| A_N\Pi^N_{\mathcal{X}} F(\bar{X})\|_{\mathcal{X}_\alpha} + \|(A - A_N) (I_{\mathcal{X}}-\Pi^N_{\mathcal{X}}) F(\bar{X})\|_{\mathcal{X}_\alpha}.
\end{align*}

Furthermore, thanks to Definition~\ref{def:A}, Lemma~\ref{lem:commute}, \eqref{hyp:B} and Lemma~\ref{lem:sum} in Appendix~\ref{appendix:A}, we are able to estimate the tail by hand as follows
\begin{align*}
    \|(A - A_N) (I_{\mathcal{X}}-\Pi^N_{\mathcal{X}}) F(\bar{X})\|_{\mathcal{X}_\alpha} 
    &=  \left\| \begin{array}{c}
    0 \\ \hline
    0 \\ \hline
    \delta (\Delta + (d-\bar{\lambda})I_{\ell^1})^{-1} (I_{\ell^1}-\Pi^N_{\ell^1}) B \bar{u} 
    \end{array}\right\|_{\mathcal{X}_\alpha} \\
    & = |\delta|\left\|(\Delta + (d-\bar{\lambda})I_{\ell^1})^{-1} (I_{\ell^1}-\Pi^N_{\ell^1}) B\bar{u} \right \|_{\alpha} \\
    & = 2|\delta| \sum_{k=N}^{+\infty} \dfrac{k^\alpha}{|-(\frac{k\pi}{l})^2 + (d-\bar{\lambda})|} \left| \sum_{j=0}^{N-1} B_{kj}\bar{u}_j \right| \\
    &\leq 
    2|\delta| \sum_{k=N}^{+\infty} \dfrac{k^\alpha}{|-(\frac{k\pi}{l})^2 + (d-\bar{\lambda})|} \sum_{j=0}^{N-1} | B_{kj}| |\bar{u}_j| \\
    &\leq 
    2C|\delta| \sum_{k=N}^{+\infty} \dfrac{k^{\alpha-q_1}}{|-(\frac{k\pi}{l})^2 + (d-\bar{\lambda})|} \sum_{j=0}^{N-1} \frac{|\bar{u}_j|}{\max(1,j^{q_2})} \\
    &\leq 2C|\delta|E(\alpha,N) \left( |\bar{u}|\cdot\chi_{q_2} \right).
\end{align*}
Thus, with $Y  =  \|A_N \Pi^N_\mathcal{X} F(\bar{X}) \|_{\mathcal{X}_\alpha} +  2C|\delta|E(\alpha,N)|\bar{u}|\cdot\chi_{q_2}$, we have $\|AF(\bar{X})\|_{\mathcal{X}_\alpha} \leq Y$.
\end{proof*}
\begin{proof*}{Derivation of $Z_1$.}
$Z_1$ is deduced from $\|I_{\mathcal{X}} -ADF(\bar{X})\|_{\mathcal{X}_\alpha}$. Let us bound $\|I_{\mathcal{X}} -ADF(\bar{X})\|_{\mathcal{X}_\alpha}$.
\begin{align*}
    I_{\mathcal{X}} -ADF(\bar{X}) = [\Pi^N_{\mathcal{X}} -\Pi^N_{\mathcal{X}}A DF(\bar{X})] + [(I_{\mathcal{X}}-\Pi^N_{\mathcal{X}}) -(I_{\mathcal{X}}-\Pi^N_{\mathcal{X}})ADF(\bar{X})].
\end{align*}
Since $\Pi^N_{\mathcal{X}}A DF(\bar{X})$ is formed by a finite number of rows and an infinite number of columns, we split again on the columns by multiplying from right by $\Pi^N_\mathcal{X}$.

Let us denote
\begin{align*}
    P_0 &= \Pi^N_{\mathcal{X}} -\Pi^N_{\mathcal{X}}A DF(\bar{X})\Pi^N_\mathcal{X} ,\\
    P_1 &= \Pi^N_{\mathcal{X}}A DF(\bar{X})(I_{\mathcal{X}} - \Pi^N_{\mathcal{X}}), \\
    P_2 &=  (I_{\mathcal{X}}-\Pi^N_{\mathcal{X}}) -(I_{\mathcal{X}}-\Pi^N_{\mathcal{X}})ADF(\bar{X}),
\end{align*}
then
\begin{align*}
    I_{\mathcal{X}} -ADF(\bar{X}) &= P_0 - P_1 + P_2.
\end{align*}

Since
\begin{align*}
    &DF(\bar{X}) =  
    \left(\begin{array}{c|c|c}
        0 & \Pi^N_{\ell^1}\tilde{u}^* & \Pi^N_{\ell^1}\tilde{v}^* \\ \hline
        -\Pi^N_{\ell^1}\bar{u} & \vartheta \Delta  + (a-\bar{\lambda})I_{\ell^1} & bI_{\ell^1} \\ \hline
        -\Pi^N_{\ell^1}\bar{v} & cI_{\ell^1} +\delta B & \Delta + (d-\bar{\lambda})I_{\ell^1}
    \end{array}\right), \text{ with } \bar{X} = \Pi^N_{\mathcal{X}}\bar{X},
\end{align*}
then, we quickly have
\begin{align*}
    P_1 &= 
    \left(\begin{array}{c|c|c}
    0 & \delta \tilde{a}^*_2B(I_{\ell^1} - \Pi^N_{\ell^1}) & 0 \\ \hline
    0 & \delta A^{12}B(I_{\ell^1} - \Pi^N_{\ell^1}) & 0 \\ \hline
    0 & \delta A^{22}B(I_{\ell^1} - \Pi^N_{\ell^1}) & 0
    \end{array} \right), \\
    P_2 &=
    \left(\begin{array}{c|c|c} 
    0 & 0 & 0 \\ \hline 
    0 & 0 & b(I_{\ell^1}-\Pi^N_{\ell^1})(\vartheta\Delta + (a-\bar{\lambda})I_{\ell^1})^{-1} \\ \hline
    0 & (I_{\ell^1} - \Pi^N_{\ell^1})(\Delta + (d-\bar{\lambda})I_{\ell^1})^{-1}(cI_{\ell^1}+\delta B) & 0
    \end{array} \right).
\end{align*}

Looking at the definition of the norm of the operator \eqref{def:op_norm}, we use the column block description.
Since $C^0_\alpha(P_1) = C^2_\alpha(P_1) = C^0_\alpha(P_2) = 0$, there are a few simplifications.

\begin{align*}
    \|P_0 - P_1 + P_2\|_{\mathcal{L}(\mathcal{X}_\alpha)} = \max \Big\{ C^0_\alpha(P_0), C^1_\alpha(P_0-P_1+P_2),  C^2_\alpha(P_0+P_2) \Big\}.
\end{align*}

Since $P_0$ is a finite matrix of size $(2N+1)\times (2N+1)$, we are able to compute its norm and its column block norms.

We then focus on $C^1_\alpha(P_0 - P_1 + P_2)$, since $(P_0-P_1+P_2)\Pi^N_\mathcal{X} = P_0 + P_2\Pi^N_\mathcal{X}$ and $(P_0-P_1+P_2)(I_\mathcal{X}-\Pi^N_\mathcal{X}) = -P_1+P_2(I-\Pi^N_\mathcal{X})$ we get,
\begin{align*}
    C^1_\alpha(P_0 - P_1 + P_2) & \leq \max \left(C^1_\alpha(P_0)+C^1_\alpha(P_2\Pi^N_\mathcal{X}), C^1_\alpha(P_1) + C^1_\alpha(P_2(I-\Pi^N_\mathcal{X}))\right).
\end{align*}

We now bound the column block norm \eqref{def:op_norm} of $P_1$, $C^1_\alpha(P_1)$,
\begin{align*}
    C^1_\alpha(P_1) 
    &= \sup_{j\geq 0}\frac{1}{\max(1,2j^\alpha)} \left\| \left(\begin{array}{c}
        \delta \tilde{a}^*_2B(I_{\ell^1} - \Pi^N_{\ell^1}) \\ \hline
        \delta A^{12}B(I_{\ell^1} - \Pi^N_{\ell^1}) \\ \hline
        \delta A^{22}B(I_{\ell^1} - \Pi^N_{\ell^1})
    \end{array} \right)_j \right\|_{\mathcal{X}_\alpha}\\
    &= |\delta| \sup_{j\geq N} \ \dfrac{1}{2j^\alpha}\left( \left|\sum_{k=0}^{N-1}\tilde{a}_{2,k}B_{kj}\right| + \left\| \left(\sum_{k=0}^{N-1} A^{12}_{ik}B_{kj} \right)_{i=0}^{N-1}  \right\|_{\ell^1_\alpha} +  \left\| \left( \sum_{k=0}^{N-1} A^{22}_{ik}B_{kj} \right)_{i=0}^{N-1}  \right\|_{\ell^1_\alpha} \right) \\
    &\leq \dfrac{C|\delta|}{2N^{\alpha+q_2}} \left( |\tilde{a}_2|\cdot\chi_{q_1} + \left\| |A^{12}|\chi_{q_1} \right\|_{\ell^1_\alpha} + \left\| |A^{22}|\chi_{q_1} \right\|_{\ell^1_\alpha} \right),
\end{align*}
we used~\eqref{hyp:B}, namely $|B_{k,j}| \leq \dfrac{C}{\max(1,k^{q_1})\max(1,j^{q_2})}$.

Similarly, we bound $C^1_\alpha(P_2\Pi^N_\mathcal{X})$ and $C^1_\alpha(P_2(I_\mathcal{X}-\Pi^N_\mathcal{X}))$.
\begin{align*}
    C^1_\alpha(P_2\Pi^N_\mathcal{X}) 
    &= \max_{j<N}\dfrac{1}{\max(1,2j^\alpha)} \left\|( (I_{\ell^1}-\Pi^N_{\ell^1})(\Delta + (d-\bar{\lambda})I_{\ell^1})^{-1}(cI_{\ell^1} + \delta B) )_j \right\|_{\ell^1_\alpha} \\
    &= \max_{j<N}\dfrac{2}{\max(1,2j^\alpha)} \sum_{k=N}^{+\infty} \left| \dfrac{c\delta_{kj}+\delta B_{kj}}{-\left(\frac{k\pi}{l}\right)^2+(d-\bar{\lambda})} \right| k^\alpha \\
    &\leq \max_{j<N} \dfrac{2C|\delta|}{\max(1,2j^{\alpha+q_2})}\sum_{k=N}^{+\infty} \dfrac{k^{\alpha-q_1}}{\left|-\left(\frac{k\pi}{l}\right)^2 + (d-\bar{\lambda})\right|}\\
    & \leq 2C|\delta| E(\alpha,N),
\end{align*}

\begin{align*}
    C^1_\alpha(P_2(I_\mathcal{X}-\Pi^N_\mathcal{X})) 
    &= \sup_{j\geq N}\dfrac{1}{\max(1,2j^\alpha)} \left\|( (I_{\ell^1}-\Pi^N_{\ell^1})(\Delta + (d-\bar{\lambda})I_{\ell^1})^{-1}(cI_{\ell^1} + \delta B) )_j \right\|_{\ell^1_\alpha} \\
    &= \sup_{j\geq N}\dfrac{2}{\max(1,2j^\alpha)} \sum_{k=N}^{+\infty} \left| \dfrac{c\delta_{kj}+\delta B_{kj}}{-\left(\frac{k\pi}{l}\right)^2+(d-\bar{\lambda})} \right| k^\alpha \\
    &\leq \sup_{j\geq N}  \dfrac{|c|}{\left|-\left(\frac{j\pi}{l}\right)^2 + (d-\bar{\lambda})\right|} +  \sup_{j\geq N} \dfrac{2C|\delta|}{\max(1,2j^{\alpha+q_2})}\sum_{k=N}^{+\infty} \dfrac{k^{\alpha-q_1}}{\left|-\left(\frac{k\pi}{l}\right)^2 + (d-\bar{\lambda})\right|}\\
    & \leq |c| R(1,d,N) + \frac{C|\delta|}{N^{\alpha+q_2}} E(\alpha,N),
\end{align*}
thanks to Lemma~\ref{lem:sum} in Appendix~\ref{appendix:A}.

Finally, we have
\begin{align}
    C^1_\alpha(P_0 - P_1 + P_2) \leq 
    & \max\Big(C^1_\alpha(P_0)+2C|\delta| E(\alpha,N), \\
    &\phantom{\max\big(} \dfrac{C|\delta|}{2N^{\alpha+q_2}} \left( |\tilde{a}_2|\cdot\chi_{q_1} + \left\| |A^{12}|\chi_{q_1} \right\|_{\ell^1_\alpha} + \left\| |A^{22}|\chi_{q_1} \right\|_{\ell^1_\alpha} \right) \\
    &\phantom{\max\big(} + |c| R(1,d,N) + \frac{C|\delta|}{N^{\alpha+q_2}} E(\alpha,N)\Big). \notag
\end{align}

Then we bound the column block norm $C^2_\alpha(P_0+P_2)$,
\begin{align}
    C^2_\alpha(P_0+P_2)  
    &\leq \max\left( C^2_\alpha(P_0)+C^2_\alpha(P_2\Pi^N_\mathcal{X}), C^2_\alpha(P_2(I_\mathcal{X}-\Pi^N_\mathcal{X}))\right) \notag \\
    &\leq \max\left(C^2_\alpha(P_0), \sup_{j\geq N}\dfrac{1}{\max(1,2j^\alpha)} \left\| (b(I_{\ell^1}-\Pi^N_{\ell^1})((\vartheta\Delta + (a-\bar{\lambda})I_{\ell^1})^{-1})_j \right\|_{\ell^1_\alpha}\right) \notag \\
    &= \max\left(C^2_\alpha(P_0), \sup_{j\geq N} \dfrac{1}{j^\alpha} \dfrac{|b| j^\alpha}{\left| - \vartheta \left(\frac{j\pi}{l}\right)^2 + (a-\bar{\lambda})\right|} \right)\notag\\
    &\leq \max\left( C^2_\alpha(P_0), |b|R(\vartheta,a,N)\right).
\end{align}

Thus,
\begin{align*}
    \|I_\mathcal{X}-ADF(\bar{X})\|_{\mathcal{L}(\mathcal{X}_\alpha)} \leq \max \{ & C^0_\alpha(P_0); \\
    & \max\big(C^1_\alpha(P_0)+2C|\delta| E(\alpha,N), \\
    &\phantom{\max\big(C^1_\alpha(P_0),}\dfrac{C|\delta|}{2N^{\alpha+q_2}} \left( |\tilde{a}_2|\cdot\chi_{q_1} + \left\| |A^{12}|\chi_{q_1} \right\|_{\ell^1_\alpha} + \left\| |A^{22}|\chi_{q_1} \right\|_{\ell^1_\alpha} \right) \\
    &\phantom{\max\big(C^1_\alpha(P_0),}+|c| R(1,d,N) + \frac{C|\delta|}{N^{\alpha+q_2}} E(\alpha,N)\big); \\ 
    & \max\left(C^2_\alpha(P_0), |b|R(\vartheta,a,N)\right) \}.
\end{align*}
Finally, since $\|\Pi^N_{\mathcal{X}} - \Pi^N_{\mathcal{X}}A DF(\bar{X})\Pi^N_\mathcal{X}\|_{\mathcal{L}(\mathcal{X}_\alpha)} = \max_{\eta\in\{ 0, 1, 2\}}\{C^\eta_\alpha(P_0)\}$ we define
\begin{align}
    Z_1 = \max\Big\{ &\|\Pi^N_{\mathcal{X}} - \Pi^N_{\mathcal{X}}A DF(\bar{X})\Pi^N_\mathcal{X}\|_{\mathcal{L}(\mathcal{X}_\alpha)} + 2C|\delta| E(\alpha,N); \\  
    &\dfrac{C|\delta|}{2N^{\alpha+q_2}} \left( |\tilde{a}_2|\cdot\chi_{q_1} + \left\| |A^{12}|\chi_{q_1} \right\|_{\ell^1_\alpha} + \left\| |A^{22}|\chi_{q_1} \right\|_{\ell^1_\alpha} \right) \notag \\
    &\quad + |c| R(1,d,N) + \frac{C|\delta|}{N^{\alpha+q_2}} E(\alpha,N) ; \notag \\
    &|b| R(\vartheta,a,N) \Big\}, \notag
\end{align}
which satisfies $\|I_{\mathcal{X}}-ADF(\bar{X})\|_{\mathcal{L}(\mathcal{X}_\alpha)} \leq Z_1$.
\end{proof*}

\begin{proof*}{Derivation of $Z_2$.}
$Z_2$ is deduced from $\|A\|_{\mathcal{L}(\mathcal{X}_\alpha)}$, we compute the operator norm of $A$ from \eqref{def:op_norm}. 
We have
\begin{align*}
    \|A\|_{\mathcal{L}(\mathcal{X}_\alpha)} &= \max_{\eta\in\{ 0, 1, 2\}} \{ C^\eta_\alpha(A)\}.
\end{align*}

Firstly, $C^0_\alpha(A) = C^0_\alpha(A_N)$ and we are able to compute it.

Secondly, for $\eta = 1,2$,
\begin{align*}
    C^\eta_\alpha(A) 
    &= C^\eta_\alpha(A_N + A-A_N)\\
    &\leq \max\left(C^\eta_\alpha(A_N), C^\eta_\alpha(A-A_N)\right) \\
    &\leq \max\left( C^\eta_\alpha(A_N),\sup_{j\in\N}\dfrac{1}{\max(1,2j^\alpha)}\left\| \left( (I_{\ell^1}-\Pi^N_{\ell^1})(\theta_\eta\Delta+(x_\eta-\bar{\lambda})I_{\ell^1})^{-1} \right)_j \right\|_{\ell^1_\alpha}\right)\\
    &\leq \max\left( C^\eta_\alpha(A_N),\sup_{j\geq N}\dfrac{1}{\max(1,2j^\alpha)}\left\| \left( (I_{\ell^1}-\Pi^N_{\ell^1})(\theta_\eta\Delta+(x_\eta-\bar{\lambda})I_{\ell^1})^{-1} \right)_j \right\|_{\ell^1_\alpha}\right)\\
    &= \max\left(C^\eta_\alpha(A_N), R(\theta_\eta,N,x_\eta)\right),
\end{align*}
with $\theta_1 = \vartheta$, $\theta_2 = 1$, $x_1 = a$, $x_2 = d$.

Thus, 
\begin{equation*}
    Z_2 = \max\{\|A_N\|_{\mathcal{L}(\mathcal{X}_\alpha)}; R(\vartheta,a,N); R(1,d,N)\},
\end{equation*}
$Z_2$ satisfies $\|A\|_{\mathcal{L}(\mathcal{X}_\alpha)} \leq Z_2$.
\end{proof*}\\
The proof of Proposition~\ref{prop:bounds} is complete.
\end{proof}

\subsubsection{Proof of Proposition~\ref{prop:lambda1}}
\label{sec:prooflambda1}

Finally we prove Proposition~\ref{prop:lambda1}, on the framing of $\delta \in [\delta_0,\delta_1] \mapsto d_0(\delta)$ by two constant piece-wise curves, by applying Theorem~\ref{thm:NK} with the explicit bounds shown in Proposition~\ref{prop:bounds}.

\begin{proof}[Proof of Proposition~\ref{prop:lambda1}]
    Let $\delta = \left[\dfrac{4k}{1000},\dfrac{4(k+1)}{1000}\right]$ with $k=607$. 
    Let $\bar{X} = (\bar{\lambda},\bar{u},\bar{v})$ be 
    the finite approximate zero of $F$ associated to $\delta$, stored in \texttt{Ubartilde(:,1)} in the file \newline \texttt{first\_eigv\_wrt\_delta\_final\_1e3.mat}, see \cite{GitHub} for technical details. We have $\bar{\lambda} = -3.1 \times 10^{-3}$. 
    Thanks to Proposition~\ref{prop:bounds}, we get that $Y = 1.301 \times 10^{-3}$, $Z_1 = 8.356\times 10^{-2}$ and $Z_2 = 1.658$ satisfy the assumptions~\eqref{hyp:NK-Y}-\eqref{hyp:NK-Z2} of Theorem~\ref{thm:NK}. Since $(1-Z_1)^2 - 4YZ_2 >0$ and $Z_1<1$, there exists $\widehat{X}\in\mathcal{B}_{\mathcal{X}_\alpha}(\bar{X},r)$ for all $r \in [r_{\min}, r_{\max}]$, $r_{\min} = 1.47 \times 10^{-3}$ and $r_{\max} = 5.53 \times 10^{-1}$. In particular, for all $\delta\in \left[\dfrac{4k}{1000},\dfrac{4(k+1)}{1000}\right]$ with $k=607$, there exists an eigenvalue $\widehat{\lambda} = \widehat{\lambda}(\delta)$ of $M$, which satisfies, $\vert \widehat{\lambda} - \bar{\lambda} \vert \leq r_{\min}$. 
    Since $\bar{\lambda} -r_{\min} > \mu $, where $\mu$ is the threshold from Theorem~\ref{thm:spectrum} depicted in Figure~\ref{fig:d0}, we have indeed enclosed the correct eigenvalue, i.e., $\widehat{\lambda} = d_0$. 
    From the values of $ \bar{\lambda}$ and $r_{\min}$ we get $-4.5\times 10^{-3} \leq d_0(\delta) \leq-1.5 \times 10^{-3}$, for all $\delta\in \left[\dfrac{4k}{1000},\dfrac{4(k+1)}{1000}\right]$ with $k=607$.
    We just showed that $\underline{d_0} : \delta\in\left[\frac{4\times607}{1000},\frac{4\times608}{1000}\right] \mapsto -4.5\times 10^{-3}$ and  $\overline{d_0} : \delta\in\left[\frac{4\times607}{1000},\frac{4\times608}{1000}\right] \mapsto-1.5\times 10^{-3}$ satisfy~\eqref{eq:enclose} on $\left[\frac{4\times607}{1000},\frac{4\times608}{1000}\right]$.

    Then, we repeat the same procedure for all $\delta = \left[\dfrac{4k}{1000},\dfrac{4(k+1)}{1000}\right]$ with $k=608, \dots, 614$, which yields the definition of the two piece-wise constant functions $\underline{d_0}$ and $\overline{d_0}$ on the entire interval $[\delta_0,\delta_1]$.

    The result is illustrated in Figure~\ref{fig:d0_focus} and given in Table~\ref{tab:enclose}. 
    The evaluation of all the bounds and the rigorous construction of $\underline{d_0}$ and $\overline{d_0}$ can be reproduced using the code available  at~\cite{GitHub}.
\end{proof}

\begin{remark}
    The values $\bar{\lambda}$, $Y$, $Z_1$, $Z_2$, $r_{\min}$, $r_{\max}$ are rounded to be easier to read. 
    The whole calculations and exhaustive results can be reproduced using the code available at~\cite{GitHub}. 
\end{remark}

\subsection{The map \texorpdfstring{$\delta \mapsto d_0(\delta)$}{}. Second properties} \label{sec:ImpTh}

In this section, we state and prove the following proposition on $\delta \mapsto d_0(\delta)$, which establishes the second part of Proposition~\ref{prop:lambda'}.

\begin{proposition} \label{prop:lambda2}
    The function $\delta \mapsto d_0(\delta)$ is piece-wise $C^1$ on $[\delta_0,\delta_1] = \bigcup_{k=607}^{614} [\frac{4k}{1000},\frac{4(k+1)}{1000}]$ and Property~\eqref{eq:lambda'} on the positivity of its derivative is satisfied.
\end{proposition}

We obtain this control on the derivative of $d_0$ thanks to the implicit function theorem. 
Indeed, the bounds obtained in Section~\ref{sec:NK} enable us to check the assumptions of the implicit function theorem, as explained in Section~\ref{sec:DFinv}, and to then construct an explicit approximation of $\delta_0'$ with computable error bounds in Section~\ref{sec:boundsder}, leading to the proof of Proposition~\ref{prop:lambda2} presented in Section~\ref{sec:prooflambda2}.

\subsubsection{Derivative of \texorpdfstring{$\delta \mapsto d_0(\delta)$}{}}
\label{sec:DFinv}

We denote by $F(\delta,X)$ the exactly same $F$ as before, but highlighting the dependency with respect to $\delta$,
\begin{equation*}
    F(\delta,X) = \begin{pmatrix}
        \begin{pmatrix} \tilde{u}_\delta \\  \tilde{v}_\delta \end{pmatrix} \cdot \begin{pmatrix} u \\  v \end{pmatrix} - 1 \\
        M(\delta)\begin{pmatrix} u \\  v \end{pmatrix} - \lambda \begin{pmatrix} u \\  v \end{pmatrix}
    \end{pmatrix}, \ X = (\lambda, u,v) \in \mathcal{X}_\alpha.
\end{equation*}
We recall that $\tilde{u}_\delta$ and $\tilde{v}_\delta$ belong to $\Pi^N_{\ell^1} \ell^1_\alpha$. 
They are piece-wise constant on $[\delta_0,\delta_1]$, since we fixed them on each $\left[\frac{4k}{1000}, \frac{4(k+1)}{1000}\right]$, $k \in \{607,\dots,614\}$. 
And $M(\delta)$ is a linear function in $\delta$.

Let us consider the function $\delta \mapsto \widehat{X}(\delta) = (d_0(\delta),\widehat{u}(\delta),\widehat{v}(\delta)) \in \mathcal{X}_\alpha$, the zero of $F$ depending on $\delta$ that belongs in the $r_{\max,\delta}$-neighborhood of $\bar{X}(\delta)$, exhibited in the proof of Proposition~\ref{prop:lambda1}.

\begin{lem}
        $\forall \delta \in [\delta_0, \delta_1], \ D_XF(\delta,\widehat{X}(\delta))$ is invertible in $\mathcal{L}(\mathcal{X}_{\alpha},\mathcal{X}_{\alpha-2})$.
\end{lem}

\begin{proof}
    For all $\delta \in [\delta_0,\delta_1]$, let $\bar{X}(\delta)$ be the approximate zero of $F(\delta,\,\cdot\,)$, that we already used in Proposition~\ref{prop:lambda1}.
    From the proof of Proposition~\ref{prop:lambda1}, we already built $A_\delta$ and we know the bounds $Y_\delta,\, Z_{1,\delta} \text{ and } Z_{2,\delta}$ satisfying the hypotheses of Theorem~\ref{thm:NK}, and yielding some $r_{\min,\delta}$. In particular, $r_{\min,\delta}$ satisfies $Z_{1,\delta} + Z_{2,\delta}r_{\min,\delta} < 1$, and therefore,
    \begin{align*}
        \|I_\mathcal{X}-A_\delta D_XF(\delta,\widehat{X}(\delta)) \|_{\mathcal{X}_\alpha} &= \|I_\mathcal{X}-A_\delta \Big(D_XF(\delta,\widehat{X}(\delta)) - D_XF(\delta,\bar{X}(\delta)) + D_XF(\delta,\bar{X}(\delta)) \Big) \|_{\mathcal{X}_\alpha}\\
        &\leq \|I_\mathcal{X}-A_\delta D_XF(\delta,\bar{X}(\delta)) \|_{\mathcal{X}_\alpha} \\
        &\quad + \|A_\delta \Big(D_XF(\delta,\widehat{X}(\delta)) - D_XF(\delta,\bar{X}(\delta)) \Big)\|_{\mathcal{X}_\alpha} \\
        &\leq Z_{1,\delta} + Z_{2,\delta}r_{\min,\delta} \\
        & < 1.
    \end{align*}
It means that $A_\delta D_XF(\delta,\widehat{X}(\delta))$ is invertible in $\mathcal{L}(\mathcal{X}_\alpha)$, and since $A_\delta$ is invertible in $\mathcal{L}(\mathcal{X}_{\alpha-2},\mathcal{X}_\alpha)$ from Definition~\ref{def:A} so is $D_XF(\delta,\widehat{X}(\delta))$ in $\mathcal{L}(\mathcal{X}_{\alpha},\mathcal{X}_{\alpha-2})$.
\end{proof}

In order to prove Proposition~\ref{prop:lambda2}, we need to study the derivative of $F$ with respect to $\delta$ on each $\left[\frac{4k}{1000},\frac{4(k+1)}{1000}\right]$ for $k \in \{ 607, \dots, 614 \}$. 

\begin{lem} \label{lem:delta_C1}
    The function $\delta \in [\delta_0,\delta_1] \mapsto \widehat{X}(\delta)$ is piece-wise $C^1$ on $\bigcup_{k=607}^{614} \left[\frac{4k}{1000}, \frac{4(k+1)}{1000}\right]$. 
    Furthermore, for all $\delta \in [\delta_0,\delta_1]$,
    \begin{equation}
        \dfrac{d\widehat{X}(\delta)}{d\delta} =-(D_{X}F(\delta,\widehat{X}(\delta)))^{-1}\left(\begin{array}{c}
        0 \\ \hline
        0 \\ \hline
        B\widehat{u}(\delta) 
    \end{array} \right).
    \end{equation}
\end{lem}

\begin{proof}
    Let $k \in \{ 607, \dots, 614 \}$.
    Let $\delta \in \left[\frac{4k}{1000},\frac{4(k+1)}{1000} \right]$, since $(\delta,\widehat{X}(\delta))$ is a zero of $F$ that is $C^1$ on $\left[ \frac{4k}{1000},\frac{4(k+1)}{1000} \right]$, and  $D_X F(\delta,\widehat{X}(\delta))$ is invertible, thanks to the implicit function theorem, we have directly that $\delta \in \left[ \frac{4k}{1000},\frac{4(k+1)}{1000} \right] \mapsto \widehat{X}(\delta)$ is $C^1$. 
    We then express $\dfrac{d\widehat{X}}{d\delta}(\delta)$.
    We first compute the derivative of $F(\delta,\widehat{X}(\delta))$ with respect to $\delta$,
    \begin{align*}
        \dfrac{dF(\delta,\widehat{X}(\delta))}{d\delta} &= \dfrac{\partial F}{\partial\delta}(\delta,\widehat{X}(\delta)) + D_{X}F(\delta,\widehat{X}(\delta))\dfrac{d\widehat{X}}{d\delta}(\delta).
    \end{align*}
    We thus have,
    \begin{align*}
        \dfrac{d\widehat{X}}{d\delta}(\delta) &=-(D_{X}F(\delta,\widehat{X}(\delta)))^{-1} \dfrac{\partial F}{\partial\delta}(\delta,\widehat{X}(\delta)),
    \end{align*}
    with
    \begin{align*}
        \dfrac{\partial F}{\partial\delta}(\delta,\widehat{X}(\delta)) &= \left(\begin{array}{c}
            0 \\ \hline
            0 \\ \hline
            B\widehat{u}(\delta) 
        \end{array} \right), \\
        D_{X}F(\delta,\widehat{X}(\delta)) &= DF(\widehat{X}(\delta)).
    \end{align*}
\end{proof}
This calculation suggests a “natural" expression of an approximation of $\dfrac{d\widehat{X}}{d\delta}(\delta)$, namely $\left[\dfrac{d\widehat{X}}{d\delta}(\delta)\right]_{app} = -A_\delta\left(\begin{array}{c}
        0 \\ \hline
        0 \\ \hline
        \Pi^N_{\ell^1}B\bar{u}(\delta)
    \end{array} \right)$, where $\bar{u}(\delta)$ comes from $\bar{X}(\delta)$, a numerical approximation of $\widehat{X}(\delta)$ (given by Theorem~\ref{thm:NK}). 
    The first element of $\left[\dfrac{d\widehat{X}}{d\delta}(\delta)\right]_{app}$ is given by a finite sum: $\left[d_0'(\delta)\right]_{app} = \tilde{a}_2^*(\delta)\Pi^N_{\ell^1}B\bar{u}(\delta)$, thanks to Definition~\ref{def:A} on $A_\delta$. 
    The element $\left[d_0'(\delta)\right]_{app}$ can therefore be computed explicitly. 
    In order to control the exact derivative of $d_0(\delta)$ with respect to $\delta$, we need to estimate

\begin{align}
    \left|d_0'(\delta) - \left[d_0'(\delta)\right]_{app}\right| 
    &\leq \left\| \dfrac{d\widehat{X}}{d\delta}(\delta) - \left[\dfrac{d\widehat{X}}{d\delta}(\delta)\right]_{app}\right\|_{\mathcal{X}_{\alpha}} \notag \\
    &\leq \left\|-D_XF(\delta,\widehat{X}(\delta))^{-1} \left(\begin{array}{c}
        0 \\ \hline
        0 \\ \hline
        B\widehat{u}(\delta)
    \end{array} \right) 
    + A_\delta\left(\begin{array}{c}
        0 \\ \hline
        0 \\ \hline
        \Pi^N_{\ell^1}B\bar{u}(\delta)
    \end{array} \right) \right\|_{\mathcal{X}_{\alpha}} .\label{eq:error_lambda_p}
\end{align}

\subsubsection{Bound on the derivative}
\label{sec:boundsder}

Here, we propose a computable bound on the left-hand side of \eqref{eq:error_lambda_p}.
We do not explicitly write the $\delta$-dependency of each object, so as not to overload the reading, especially in calculations.

\begin{lem} \label{lem:bound_lambda} 
Let $\delta \in [\delta_0, \delta_1]$.
With the notations introduced in Section~\ref{sec:NK} and Section~\ref{sec:ImpTh} we have,
    \begin{align*}
    \left|d_0' - \left[d_0'\right]_{app}\right| &\leq \dfrac{1}{1-(Z_1+Z_2 r_{\min})}\Bigg( C \left(|\tilde{a}_2|\cdot\chi_{q_1} + \|| A^{12}|\chi_{q_1}\|_{\ell^1_\alpha} + \||A^{22}|\chi_{q_1}\|_{\ell^1_\alpha} + 2E(\alpha,N) \right)r_{\min} \\
    &\qquad +  2CE(\alpha,N)|\bar{u}| \cdot \chi_{q_2} + (Z_1+Z_2r_{\min}) \left\|\left(\begin{array}{c}
            \tilde{a}_2^*\Pi^N_{\ell^1}B\bar{u} \\ \hline
            A^{12}\Pi^N_{\ell^1}B\bar{u} \\ \hline
            A^{22}\Pi^N_{\ell^1}B\bar{u}
        \end{array} \right)  
        \right\|_{\mathcal{X}_{\alpha}} \Bigg).
\end{align*}
\end{lem}

\begin{proof}
    Let $\delta \in [\delta_0,\delta_1]$.
    Firstly, by the Neumann series and Theorem~\ref{thm:NK}, we can express $D_XF(\delta,\widehat{X})^{-1}$ with objects we control. 
    We already know that $Z_1 + Z_2r_{\min}<1$, thus
    \begin{align*}
        D_XF(\delta,\widehat{X})^{-1} &= D_XF(\delta,\widehat{X})^{-1}A^{-1}A \\
        &=  (AD_XF(\delta,\widehat{X}))^{-1}A \\
        &= \sum_{k=0}^{+\infty} [I_{\mathcal{X}} - AD_XF(\delta,\widehat{X})]^k A \\
        &= \sum_{k=0}^{+\infty} \left[I_{\mathcal{X}} - AD_XF(\delta,\bar{X}) - A\big(D_XF(\delta,\widehat{X})-D_XF(\delta,\bar{X})\big)\right]^k A. \\
    \intertext{Likewise,}
        D_XF(\delta,\widehat{X})^{-1} - A &= \sum_{k=1}^{+\infty} \left[I_{\mathcal{X}} - AD_XF(\delta,\bar{X}) - A\big(D_XF(\delta,\widehat{X})-D_XF(\delta,\bar{X})\big)\right]^k A,
    \end{align*}
    which is indeed convergent since 
    \begin{align*}
        \|I_{\mathcal{X}} - AD_XF(\delta,\bar{X}) - A\big(D_XF(\delta,\widehat{X})-D_XF(\delta,\bar{X})\big)\|_{\mathcal{L}(\mathcal{X}_\alpha)} \leq Z_1 + Z_2 r_{\min}.
    \end{align*}
    Thus, for any $Z\in\mathcal{X}_\alpha$,
    \begin{align*}
        \|-(D_XF(\delta,\widehat{X}))^{-1}Z\|_{\mathcal{X}_\alpha} &\leq \dfrac{1}{1-(Z_1+Z_2 r_{\min})} \|AZ\|_{\mathcal{X}_\alpha},  
    \end{align*}
    and
    \begin{align*}
        \|(-(D_XF(\delta,\widehat{X}))^{-1} + A)Z\|_{\mathcal{X}_\alpha} &\leq \dfrac{Z_1+Z_2 r_{\min}}{1-(Z_1+Z_2 r_{\min})} \|AZ\|_{\mathcal{X}_\alpha}.   
    \end{align*}
    
    Secondly, $B\widehat{u} = B(\widehat{u}-\bar{u}) + (I-\Pi^N)B\bar{u} + \Pi^NB\bar{u}$. 
    Thus, with the triangle inequality we continue to bound \eqref{eq:error_lambda_p},
    \begin{align}
        \left|d_0' - \left[d_0'\right]_{app}\right| &\leq
        \left\|D_XF(\delta,\widehat{X})^{-1} \left(\begin{array}{c}
            0 \\ \hline
            0 \\ \hline
            B(\widehat{u}-\bar{u})
        \end{array} \right)  \right\|_{\mathcal{X}_{\alpha}} +
        \left\|D_XF(\delta,\widehat{X})^{-1} \left(\begin{array}{c}
            0 \\ \hline
            0 \\ \hline
            (I_{\ell^1}-\Pi^N_{\ell^1})B\bar{u}
        \end{array} \right) \right\|_{\mathcal{X}_{\alpha}} \notag \\
        &\hspace{3em} + \left\|(D_XF(\delta,\widehat{X})^{-1}-A)\left(\begin{array}{c}
            0 \\ \hline
            0 \\ \hline
            \Pi^N_{\ell^1}B\bar{u}
        \end{array} \right) \right\|_{\mathcal{X}_{\alpha}} \notag \\ 
        &\leq 
        \dfrac{1}{1-(Z_1+Z_2 r_{\min})}\Bigg(\left\|A\left(\begin{array}{c}
            0 \\ \hline
            0 \\ \hline
            B(\widehat{u}-\bar{u})
        \end{array} \right)  \right\|_{\mathcal{X}_{\alpha}} + \left\|A\left(\begin{array}{c}
            0 \\ \hline
            0 \\ \hline
            (I_{\ell^1}-\Pi^N_{\ell^1})B\bar{u}
        \end{array} \right)  \right\|_{\mathcal{X}_{\alpha}} \notag \\
        &\hspace{3em} +(Z_1+Z_2 r_{\min}) \left\|A\left(\begin{array}{c}
            0 \\ \hline
            0 \\ \hline
            \Pi^N_{\ell^1}B\bar{u}
        \end{array} \right)  \right\|_{\mathcal{X}_{\alpha}} \Bigg). \label{eq:error_lambda_p2}
    \end{align}
    
    Then, we bound each term of this sum. 
    The calculations are similar to those made in the proof of Proposition~\ref{prop:bounds}, so we omit details. Since $|B_{k,j}| \leq C(\chi_{q_1})_k(\chi_{q_2})_j$, $(I-\Pi^N)A$ is diagonal and $\tilde{a}_2,\ \bar{u}$ are finite, it is straightforward to get 
    \begin{align*}
        \left\|A\left(\begin{array}{c}
            0 \\ \hline
            0 \\ \hline
            B(\widehat{u}-\bar{u})
        \end{array} \right)  \right\|_{\mathcal{X}_{\alpha}}\hspace{-1em} &= \left\|\begin{array}{c}
            \tilde{a}_2^*B(\widehat{u}-\bar{u}) \\ \hline
            A^{12}B(\widehat{u}-\bar{u}) \\ \hline
            (A^{22}+(I_{\ell^1} - \Pi^N_{\ell^1})(\Delta+(d-\bar{d_0})I_{\ell^1})^{-1})B(\widehat{u}-\bar{u})
        \end{array} \right\|_{\mathcal{X}_{\alpha}} \\
        &\leq \left(C\left(|\tilde{a}_2|\cdot\chi_{q_1} + \||A^{12}|\chi_{q_1}\|_{\ell^1_\alpha} + \||A^{22}|\chi_{q_1}\|_{\ell^1_\alpha}\right) + 2CE(\alpha,N) \right)r_{\min},
    \end{align*}
    as well as
    \begin{align*}
        \left\|A\left(\begin{array}{c}
            0 \\ \hline
            0 \\ \hline
            (I_{\ell^1}-\Pi^N_{\ell^1})B\bar{u}
        \end{array} \right)  \right\|_{\mathcal{X}_{\alpha}}\hspace{-1em}
        &= \left\|\begin{array}{c}
            0 \\ \hline
            0 \\ \hline
            (I_{\ell^1} - \Pi^N_{\ell^1})(\Delta+(d-\bar{d_0})I_{\ell^1})^{-1} B\bar{u}
        \end{array} \right\|_{\mathcal{X}_{\alpha}}\hspace{-1em} \\
        &\leq 2CE(\alpha,N)|\bar{u}| \cdot \chi_{q_2},
    \end{align*}
    and
     \begin{equation*}
        \left\|A\left(\begin{array}{c}
            0 \\ \hline
            0 \\ \hline
            \Pi^N_{\ell^1}B\bar{u}
        \end{array} \right)  \right\|_{\mathcal{X}_{\alpha}}\hspace{-1em} =  \left\|\left(\begin{array}{c}
            \tilde{a}_2^*\Pi^N_{\ell^1}B\bar{u} \\ \hline
            A^{12}\Pi^N_{\ell^1}B\bar{u} \\ \hline
            A^{22}\Pi^N_{\ell^1}B\bar{u}
        \end{array} \right)  
        \right\|_{\mathcal{X}_{\alpha}}. \qedhere 
    \end{equation*}
\end{proof}

\subsubsection{Proof of Proposition~\ref{prop:lambda2}}
\label{sec:prooflambda2}
We can now complete the proof of the result announced in Section~\ref{sec:ImpTh}.

\begin{proof}[Proof of Proposition~\ref{prop:lambda2}.]
    From Lemma~\ref{lem:delta_C1}, $\delta \mapsto d_0(\delta)$ is piece-wise $C^1$ on $[\delta_0,\delta_1]$. We now establish Property~\eqref{eq:lambda'}.

    Let $\delta = \left[\dfrac{4k}{1000},\dfrac{4(k+1)}{1000}\right]$ with $k=607$. Let $\bar{X}_\delta = (\bar{\lambda}_\delta,\bar{u}_\delta,\bar{v}_\delta)$ a finite approximate zero of $F$ associated to $\delta$ that we already computed in Proposition~\ref{prop:lambda1}.
    Thanks to Proposition~\ref{prop:bounds}, we compute $Y_\delta$, $Z_{1,\delta}$ and $Z_{2,\delta}$.
    $Y_\delta = 1.301 \times 10^{-3}$, $Z_{1,\delta} = 8.356\times 10^{-2}$ and $Z_{2,\delta} = 1.658$.
    Let $\widehat{X}(\delta)$, the zero of $F(\delta,\cdot)$, such that $\widehat{X}(\delta) \in \mathcal{B}_{\mathcal{X}_\alpha}(\bar{X}_\delta,r)$ for all $r \in [r_{\min,\delta}, r_{\max,\delta}]$, with $r_{\min,\delta} = 1.47 \times 10^{-3}$ and $r_{\max,\delta} = 5.53 \times 10^{-1}$. This is obtained from Theorem~\ref{thm:NK} in the proof of Proposition~\ref{prop:lambda1}.
    From Lemma~\ref{lem:bound_lambda}, we compute the approximate value and the associated bound.
    We have 
    \begin{align*}
        \left[d_0'(\delta)\right]_{app} &= 0.21, \\ 
        \left|d_0'(\delta)-\left[d_0'(\delta)\right]_{app}\right| &\leq 0.05.
    \end{align*}
    Therefore, we have $d_0'(\delta_0) \subset [0.16, 0.26]$.
    
    Then, we repeat the same procedure for each $\delta = \left[\dfrac{4k}{1000},\dfrac{4(k+1)}{1000}\right]$ with $k=608, \dots, 614$. We obtain similar results.
    All computations and figures can be reproduced, see \cite{GitHub}.
    
    We conclude that $d_0' \left( \left[\delta_0,\delta_1 \right]\right) \subset [0.16,0.26]$.
\end{proof}

\subsection{Proof of the Existence of a threshold}

We conclude on the existence and uniqueness of a threshold $\delta^*$ in $[0,4]$.
We recap what we know.
From Section~\ref{sec:Gersh}, we have for each $\delta\in[0,4]$, the first eigenvalue (largest real part) $d_0(\delta)$ of $M$ is isolated.
We have also that, for $\delta \in [0,\delta_0]$, $d_0(\delta) < 0$, and for $\delta \in [\delta_1,4]$, $d_0(\delta) > 0$. 
And finally, we prove below Theorem~\ref{thm:threshold}.
\begin{proof}[Proof of Theorem~\ref{thm:threshold}.]
    From Proposition~\ref{prop:lambda'} we obtain $d_0'(\delta) > 0$ for all $\delta \in [\delta_0,\delta_1]$, it means $d_0$ is continuous and strictly increasing in that range, from $-1.5\times 10^{-3}$ to $1.3 \times 10^{-3}$ in the worst case, Proposition~\ref{prop:lambda1}.
    According to the mean value theorem, there exists a unique $\delta^* \in (\delta_0,\delta_1)$, such that $d_0(\delta^*) = 0$ with the conclusion announced.
\end{proof}

\section{Outlook}

The study of system~\eqref{eq:system} was motivated by a mathematical model for the dynamics of virus and $\mathsf{T}$~cells during inflammatory processes, \cite{reisch_chemotactic_2019}. Here, we proved exemplary the destabilization of a trivial steady state by nonlocal effects. In contrast, the situation in modeling liver inflammation differs from this example: Usually, the inflammatory process develops towards a chronic state that shows a spatially heterogeneous stationary spread of virus and $\mathsf{T}$~cells. Starting from this nontrivial steady state, medical treatments may change the parameter $\delta$, which is interpreted as the strength of the immune system. The objective of the therapy is to change the inflammation from the chronic state towards a state without any virus. Mathematically, this means destabilizing a stable nontrivial steady state and gaining dynamics towards the trivial steady state. 
Consequently, a future step in the research is to start the investigations from a nontrivial steady state. 

Compared to the current work, one would first need to obtain a rigorous and accurate description of such a nontrivial steady state. To that end, a computer-assisted argument based on the Newton-Kantorovich Theorem~\ref{thm:NK} could prove suitable, with $F$ describing the stationary problem. One should then be able to use a Gershgorin argument as in Section~\ref{sec:Gersh} to precisely study the spectrum of the linearization at that nontrivial steady state. However, even if one can prove that a nontrivial steady state destabilizes when $\delta$ is appropriately varied, whether the system then converges back to the trivial equilibrium or not is a delicate question, which may be attacked using the rigorous integrators recently developed in~\cite{BerBreShe24,WilZgl24,DucLesTak25}. These questions will be the subject of further investigations.

\appendix
\section{Appendix} \label{appendix:A}
We go back to the system \eqref{eq:system_infinite}, with $B_{i,j} =  \fint_{\Omega_1} \varphi_i(x)\mathrm{d}x\int_{\Omega_2} \varphi_j (x)\mathrm{d}x$.
The next proposition establishes hypotheses~\eqref{hyp:B} for this $B$.

\begin{proposition} \label{prop:Bcoef}
Let $\Omega = [x_{\min}, x_{\max} ], \Omega_1,\Omega_2 \subset \Omega$ such that 
$\Omega_1 = \bigsqcup_{k=1}^{I_1} [a^k_1,b^k_1]$, $\Omega_2 = \bigsqcup_{k=1}^{I_2} [a^k_2,b^k_2]$ (disjoint union).
We have for all $i,j\in \N$, $|B_{i,j}|\leq \dfrac{C}{\max(1,i)\max(1,j)}$,

where $C = \max \left( \dfrac{8I_1I_2|\Omega|}{\pi^2 |\Omega_1|}, \dfrac{4I_2}{\pi}, \dfrac{2I_1|\Omega_2|}{|\Omega_1|\pi}, \dfrac{|\Omega_2|}{|\Omega|} \right).$
\end{proposition}
\begin{proof}
    Firstly, we will study the linear operator $B : u(\cdot) \longmapsto \left(\dfrac{1}{|\Omega_1|}\int_{\Omega_2} u(x) dx \right) \1_{\Omega_1}(\cdot)$  through the Fourier modes.
    Let $j\in\N$, $u_j = \displaystyle \dfrac{1}{|\Omega|} \int_{|\Omega|} u(x) \cos\left(\dfrac{x-x_{\min}}{x_{\max}-x_{\min}}\pi \cdot j \right)dx$, the $j^{th}$ Fourier mode of $u$.
    We call $\mathcal{F}(u) = \left( u_j \right)_{j\in\N}$, the sequence of Fourier's coefficient. 
    We have 
    \begin{align*}
        \displaystyle u(x) = u_0 + 2\sum_{j=1}^{+\infty} u_j \cos\left(\dfrac{x-x_{\min}}{x_{\max}-x_{\min}}\pi \cdot j \right).
    \end{align*}

    Let $i\in\N$,
    \begin{align*}
        \mathcal{F}_i \left(B(u)\right) &= \left(\dfrac{1}{|\Omega_1|}\int_{\Omega_2} u(x) dx \right) \mathcal{F}_i(\1_{\Omega_1}) \\
        & = \dfrac{1}{|\Omega_1|} \left( u_0 |\Omega_2| + 2\sum_{j=1}^{+\infty} u_j \int_{\Omega_2} \cos \left( \dfrac{x-x_{\min}}{x_{\max}-x_{\min}}\pi \cdot j \right) dx \right) \mathcal{F}_i(\1_{\Omega_1}) \\
        &= \dfrac{1}{|\Omega_1|} |\Omega| \left( \widetilde{\mathcal{F}(\1_{\Omega_2})} \cdot \mathcal{F}(u) \right) \mathcal{F}_i(\1_{\Omega_1}),
    \end{align*}
    where $\widetilde{\mathcal{F}(\1_{\Omega_2})}$ is the sequence defined as $\forall k \in \N, \ \widetilde{\mathcal{F}(\1_{\Omega_2})}_k = \left \{ \begin{array}{ll}
    \mathcal{F}_0(\1_{\Omega_2}) & k=0   \\
     2\mathcal{F}_k(\1_{\Omega_2}) & k \geq 1.
    \end{array}\right.$

    We integrate the cosine functions on $\Omega_1$ to get $\forall i\in\N$,
    \begin{align*}
        \mathcal{F}_i(\1_{\Omega_1}) &= \left\{ \begin{array}{l}
            \dfrac{|\Omega_1|}{|\Omega|}, \ i=0 \\
            \displaystyle \dfrac{1}{\pi\cdot i} \sum_{k=1}^{I_1}\sin\left(\dfrac{b_1^k - x_{\min}}{x_{\max}-x_{\min}}\pi \cdot i\right) - \sin\left(\dfrac{a_1^k - x_{\min}}{x_{\max}-x_{\min}}\pi \cdot i\right), \   i\geq 1.  
        \end{array} \right.
    \end{align*}
    We can write a similar result for $\mathcal{F}(\1_{\Omega_2})$.

Finally, let $i,j\in\N$, the coefficient $B_{i,j}$ is $B$ applied to the $j^{th}$ mode, $x \mapsto \cos\left(\frac{x-x_{min}}{x_{max}-x_{min}}\pi\cdot j\right)$, projected in the $i^{th}$ mode. 
We have
    \begin{align*}
        B_{i,j} &= \dfrac{|\Omega|}{|\Omega_1|} \mathcal{F}_i(\1_{\Omega_1})\times \widetilde{\mathcal{F}_j(\1_{\Omega_2})}.
    \end{align*}
We verify that $ (B_{i,j})_{j\in\N}\cdot \mathcal{F}(u) = \mathcal{F}_i(B(u))$.

Then, for $i \geq 1$ and $j \geq 1$, we have
    \begin{align}
        B_{i,j} =\dfrac{|\Omega|}{|\Omega_1|} &\left(\dfrac{1}{\pi \cdot i} \sum_{k=1}^{I_1}\sin\left(\dfrac{b_1^k - x_{\min}}{x_{\max}-x_{\min}}\pi \cdot i\right) - \sin\left(\dfrac{a_1^k - x_{\min}}{x_{\max}-x_{\min}}\pi \cdot i\right) \right) \notag\\ 
        \times &\left( \dfrac{2}{\pi \cdot j} \sum_{k=1}^{I_2}\sin\left(\dfrac{b_2^k - x_{\min}}{x_{\max}-x_{\min}}\pi \cdot j\right) - \sin\left(\dfrac{a_2^k - x_{\min}}{x_{\max}-x_{\min}}\pi \cdot j\right) \right). \label{eq:Bij}
    \end{align}
Then,
\begin{align*}
        |B_{i,j}| \leq \dfrac{|\Omega|}{|\Omega_1|} \dfrac{2I_1}{\pi \cdot i} \dfrac{4I_2}{\pi \cdot j}
\end{align*}

When $i=0$, we replace in \eqref{eq:Bij}, the first parenthesis by $\dfrac{|\Omega_1|}{|\Omega|}$. 
We have $$|B_{0j}| \leq \dfrac{|\Omega|}{|\Omega_1|} \dfrac{|\Omega_1|}{|\Omega|} \dfrac{4I_2}{\pi \cdot j} = \dfrac{4I_2}{\pi \cdot j}. $$

When $j=0$, we replace in \eqref{eq:Bij} the second parenthesis by $\dfrac{|\Omega_2|}{|\Omega|}$. 
We have $$|B_{i0}| \leq \dfrac{|\Omega|}{|\Omega_1|} \dfrac{2I_1}{\pi \cdot i} \dfrac{|\Omega_2|}{|\Omega|} = \dfrac{2I_1|\Omega_2|}{|\Omega_1| \pi \cdot i}. $$

Finally, for $i=0,\ j=0$, we have $$|B_{00}| \leq \dfrac{|\Omega|}{|\Omega_1|} \dfrac{|\Omega_1|}{|\Omega|} \dfrac{|\Omega_2|}{|\Omega|} = \dfrac{|\Omega_2|}{|\Omega|}.$$

Combining all, we have the result for all $i,j \in \N$.
\end{proof}

\begin{lem} \label{lem:inequality}
Let $f : x \mapsto \max(1,x^p)$, with $p>1-q$.
    $$\displaystyle \sum_{j=N}^{+\infty} \frac{1}{f(j)j^q} \leq \dfrac{1}{(p+q-1)(N-1)^{p+q-1}}$$ 
\end{lem}
\begin{proof}
The function $x \mapsto \dfrac{1}{x^{p+q}}$ is positive and decreasing, so

$\forall j\in \N, \ j\geq N, \ \forall x\in[j,j+1]  \quad \dfrac{1}{(j+1)^{p+q}}\leq \dfrac{1}{x^{p+q}}\leq \dfrac{1}{j^{p+q}}$.

We deduce that
\begin{align*}
    \int_{N}^{+\infty} \dfrac{1}{x^{p+q}} dx &\leq \sum_{j=N}^{+\infty} \frac{1}{j^{p+q}}
    \leq \int_{N-1}^{+\infty} \dfrac{1}{x^{p+q}} dx.
\end{align*}
So,
\begin{align*}
    \int_{N}^{+\infty}x^{-(p+q)} dx &\leq \sum_{j=N}^{+\infty} \frac{1}{j^{p+q}} 
    \leq \int_{N-1}^{+\infty} x^{-(p+q)} dx.
\end{align*}
So,
\begin{align*}
    \dfrac{1}{(p+q-1)N^{p+q-1}} &\leq \sum_{j=N}^{+\infty} \frac{1}{j^{p+q}} 
    \leq \dfrac{1}{(p+q-1)(N-1)^{p+q-1}}.
\end{align*}
\end{proof}

\begin{lem} \label{lem:sum}
Let $\alpha \in [0,1+q), \ N\in\N,\ N > 1+\frac{l}{\pi}\sqrt{|d-\lambda|}$. Let $E(\alpha,N)$ defined in Proposition~\ref{prop:bounds}. 
We have,
\begin{align*}
     \sum_{k=N}^{+\infty} \dfrac{k^{\alpha-q}}{\left|-\left(\dfrac{k\pi}{l}\right)^2 + (d-\lambda)\right|} &\leq E(\alpha,N)
\end{align*}
\end{lem}
\begin{proof}
Firstly, assume $d-\lambda \leq 0$, let $k \geq N$, we have $\displaystyle\left|-\left(\frac{k\pi}{l}\right)^2 + (d-\lambda)\right| = \left(\frac{k\pi}{l}\right)^2 + (\lambda-d) \geq \left(\frac{k\pi}{l}\right)^2$.
Thus, $\displaystyle\sum_{k=N}^{+\infty} \dfrac{k^{\alpha-q}}{|-(\frac{k\pi}{l})^2 + (d-\lambda)|} \leq \left(\frac{l}{\pi}\right)^2 \sum_{k=N}^{+\infty} k^{\alpha-q-2}$, which converges.
And by integral comparison,

$\displaystyle\sum_{k=N}^{+\infty} \dfrac{k^{\alpha-q}}{|-(\frac{k\pi}{l})^2 + (d-\lambda)|} \leq \left(\frac{l}{\pi}\right)^2  \int_{N-1}^{+\infty} x^{\alpha-q-2}dx = \left(\frac{l}{\pi}\right)^2 \dfrac{(N-1)^{\alpha-(1+q)}}{(1+q)-\alpha} = E(\alpha,N)$, since $\nu = 0$.

Then, assume  $d-\lambda > 0$, we have 
\begin{align*}
    \left|-\left(\dfrac{k\pi}{l}\right)^2 + (d-\lambda)\right| &= \left|\left(\dfrac{k\pi}{l}\right)^2 - (d-\lambda)\right| \\
    &= \left(\dfrac{\pi}{l}\right)^2 \left|k-\dfrac{l}{\pi}\sqrt{d-\lambda}\right| \left|k+\dfrac{l}{\pi}\sqrt{d-\lambda}\right| \\
    &\geq  \left(\dfrac{\pi}{l}\right)^2 \left|k-\dfrac{l}{\pi}\sqrt{d-\lambda}\right| k .
\end{align*} 
In this case, recall that $\nu = \frac{l}{\pi}\sqrt{d-\lambda}$,  we have $N\geq\nu$ then $\left|-\left(\frac{k\pi}{l}\right)^2 + (d-\lambda)\right| =   \left(\frac{\pi}{l}\right)^2 (k-\nu) k$.
Thus, $\displaystyle \dfrac{k^{\alpha-q}}{|-\left(\frac{k\pi}{l}\right)^2 + (d-\lambda)|} \leq  \left(\frac{l}{\pi}\right)^2\dfrac{k^{\alpha-q-1}}{k-\nu}$.

By integral comparison,
\begin{align*}
    \sum_{k=N}^{+\infty} \dfrac{k^{\alpha-q}}{|-(\frac{k\pi}{l})^2 + (d-\lambda)|} &\leq \left(\frac{l}{\pi}\right)^2  \int_{N-1}^{+\infty} \dfrac{x^{\alpha-q-1}}{x-\nu}dx \\
    &\leq \left(\frac{l}{\pi}\right)^2  \int_{N-1-\nu}^{+\infty} \dfrac{(y+\nu)^{\alpha-q-1}}{y}dy \\
    \intertext{since $\alpha < 1+q$, it means $y \mapsto y^{\alpha-q-1}$ decreases, and $\nu>0$ we have}
    &\leq \left(\frac{l}{\pi}\right)^2  \int_{N-1-\nu}^{+\infty} y^{\alpha-q-2}dy \\
    &= \left(\frac{l}{\pi}\right)^2 \dfrac{(N-1-\nu)^{\alpha-(1+q)}}{(1+q)-\alpha}\\
    &= E(\alpha,N).
\end{align*}
\end{proof}

\section*{Acknowledgment}

This project started during a visit of MB and MP at the University of Graz, supported by the French-Austrian Amadeus project \emph{Fractional cross-diffusion equations}. MB and MP also acknowledges the support of the ANR project CAPPS: ANR-23-CE40-0004-01.
CR acknowledges the support by the German Academic Exchange Service (DAAD) through a fellowship at the University of Graz.  BT is supported by the FWF project ``Quasi-steady-state approximation for PDE''
number I-5213. This work is partially supported by NAWI Graz.

\bibliographystyle{alpha}
\bibliography{Ref}
\end{document}